\documentclass{amsart} 

\usepackage[english]{babel}
\usepackage[latin1]{inputenc} 
\usepackage{amsmath}
\usepackage{amsthm}
\usepackage{amssymb}
\usepackage{tikz}

\usepackage{imakeidx}
\usepackage{appendix}
\usepackage{tikz-cd}
\usepackage{verbatim}
\usepackage{enumitem}
\usepackage{multicol}

\usepackage{hyperref}
\usepackage{cleveref}
\usepackage{xr}

\newtheorem{thm}{Theorem}[section] 
\newtheorem{prop}[thm]{Proposition}

\newtheorem{lemma}[thm]{Lemma}
\newtheorem{cor}[thm]{Corollary}

\newtheorem{fact}[thm]{Fact}
\newtheorem{conj}[thm]{Conjecture}

\theoremstyle{definition}

\newtheorem{example}[thm]{Example}
\newtheorem{claim}{Claim}
\theoremstyle{remark}

\newtheorem{rmk}[thm]{Remark}

\numberwithin{equation}{section} 

\newcommand{\R}{\mathbb R}

\newcommand{\Q}{\mathbb Q}

\newcommand{\G}{\mathbb G}
\newcommand{\Z}{\mathbb Z}

\renewcommand{\P}{\mathbb P}

\renewcommand{\c}{\subseteq}

\renewcommand{\O}{\mathcal O}

\newcommand{\N}{\mathbb N}
\newcommand{\A}{\mathbb A}

\newcommand{\mc}[1]{\mathcal{#1}}
\newcommand{\cl}{\overline}
\newcommand{\set}[1]{\{#1\}}

\renewcommand{\phi}{\varphi}
\newcommand{\on}[1]{\operatorname{#1}}

\DeclareMathOperator{\Spec}{Spec}

\DeclareMathOperator{\im}{im}

\DeclareMathOperator{\Gm}{\mathbb{G}_m}

\DeclareMathOperator{\ed}{ed}
\DeclareMathOperator{\ind}{ind}
\DeclareMathOperator{\trdeg}{trdeg}

\newcommand{\ang}[1]{\left \langle{#1}\right \rangle}
\newcommand{\floor}[1]{\left \lfloor{#1}\right \rfloor}

\makeatletter
\newcommand{\doublewidetilde}[1]{{%
		\mathpalette\double@widetilde{#1}%
}}
\newcommand{\double@widetilde}[2]{%
	\sbox\z@{$\m@th#1\widetilde{#2}$}%
	\ht\z@=.9\ht\z@
	\widetilde{\box\z@}%
}
\makeatother

\title{Essential dimension and genericity for quiver representations}
\author{Federico Scavia}

\begin{document}
	
	\begin{abstract}
		We study the essential dimension of representations of a fixed quiver with given dimension vector. We also consider the question of when the genericity property holds, i.e., when essential dimension and generic essential dimension agree. We classify
		the quivers satisfying the genericity property for every dimension
		vector and show that for every wild quiver the genericity property holds
		for infinitely many of its Schur roots. We also construct a large class
		of examples, where the genericity property fails. Our results are
		particularly detailed in the case of Kronecker quivers.
	\end{abstract}
	\maketitle
	
	\section{Introduction}
	Given an algebra $A$, it is a natural goal to understand the category of its representations, and if possible to give a classification. Initially one would like to describe representations over an algebraically closed field. However, it is also interesting to study representations of $A$ over non-algebraically closed fields. A template for this approach is provided by the classical theory of representations of finite groups (or equivalently, their group algebras), as summarized, e.g., in the books \cite{curtis1966representation} or \cite{serre2012linear}. In particular, it is interesting to understand which representations are defined over which fields. This leads to the study of essential dimension in representation theory; see \cite{karpenko2014numerical}, \cite{bensonreichstein} and \cite{scaviaalgebras}. 
	
	In this paper we will focus on representations of quiver path algebras. This is a large and interesting family of algebras, which has found numerous applications in algebraic geometry, Lie theory and physics. An important distinguishing feature of this family of algebras is that here representation-theoretic results can often be expressed in combinatorial (graph-theoretic) language. We initiated the study of essential dimension of quiver representations in the second half of \cite{scaviaalgebras}. This
	paper is a sequel to \cite{scaviaalgebras}, with a focus on the genericity property.
	
	Let $k$ be a field. Following P. Brosnan, Z. Reichstein and A. Vistoli \cite{genericity0}, we define the essential dimension of an algebraic stack $\mc{X}$ over $k$ as the minimal number $\ed_k\mc{X}$ of parameters required to describe any object of $\mc{X}$. If $\mc{X}$ is integral, we define the generic essential dimension $\on{ged}_k\mc{X}$ as the essential dimension of a generic object of $\mc{X}$. We say that the genericity property holds for $\mc{X}$ if $\on{ged}_k\mc{X}=\ed_k\mc{X}$; see \Cref{essdim} for the precise definitions.
	
	The genericity property fails in general (see \cite[Example 6.5]{genericity0}) but holds for smooth algebraic stacks with reductive automorphism groups \cite{genericity} (and in particular, Deligne-Mumford stacks \cite{genericity0}). In many interesting examples where these conditions are not satisfied, the genericity property continues to hold \cite{biswasdhillonhoffmann, genericity}. This phenomenon is poorly understood; one of the goals of this paper is to investigate the genericity property of stacks of quiver representations. In particular, we produce large families of examples where genericity holds and where it fails.
	
	Representations of dimension $\alpha$ of a fixed quiver $Q$ are parametrized by an integral stack $\mc{R}_{Q,\alpha}$ of finite type over $k$ (see \Cref{essdim}), and it makes sense to consider its generic essential dimension. In \Cref{adhoc} we give an equivalent definition of $\on{ged}_k\mc{R}_{Q,\alpha}$, not involving stacks.
	
	In this work, we study $\on{ged}_k\mc{R}_{Q,\alpha}$ and the genericity property for $\mc{R}_{Q,\alpha}$. On the one hand, this improves our understanding of the essential dimension of representations of algebras. On the other hand, this is the first appearance of a large family of counterexamples to the genericity property. The algebraic stacks $\mc{R}_{Q,\alpha}$ are smooth, but their automorphism groups are often non-reductive, and so it is natural to investigate what happens in this case.
	
	Our first result summarizes our understanding of the generic essential dimension of $\mc{R}_{Q,\alpha}$. We refer the reader to \Cref{quiverreps} for the necessary definitions of the Tits form and Schur roots.

	\begin{thm}\label{genericed}
		Let $Q$ be a quiver, and let $\alpha$ be a Schur root of $Q$. 
		We have \begin{equation}\label{coll1}\on{ged}_k\mc{R}_{Q,\alpha}\leq 1-\ang{\alpha,\alpha}+\sum_{p}(p^{v_p(\gcd(\alpha_i))}-1),\end{equation} where the sum is over all prime numbers $p$. One has equality if \Cref{collconj} holds for $d=\gcd(\alpha_i)$.
	\end{thm}
	We generalize this result to the case when $\alpha$ is an arbitrary root of $Q$ in \Cref{arbitraryroot}.

	When the genericity property holds, the same formulas are true for the essential dimension. It is then natural to wonder when the genericity property holds. We have two results in this direction. 
	
	\begin{thm}\label{genericityproperty}
		Let $Q$ be a quiver. Then $\mc{R}_{Q,\alpha}$ satisfies the genericity property for every dimension vector $\alpha$ if and only if $Q$ is of finite representation type or has at least one loop at every vertex.
	\end{thm}
	As an important special case, the combination of \Cref{genericed} and \Cref{genericityproperty} gives us a formula for the essential dimension of the $n$-dimensional representations of the $r$-loop quiver; see \Cref{exloop}. 
	
	\begin{thm}\label{infiniteged}
		Let $Q$ be a wild quiver. There are infinitely many Schur roots $\alpha$ such that the genericity property holds for $\mc{R}_{Q,\alpha}$.
	\end{thm}
	For a constructive variant of this result, see \Cref{moreprecise}.
	
	Our final result concerns generalized Kronecker quivers:
	\[
	\begin{tikzcd}[row sep=1ex]
	& 1 \arrow[r] \arrow[r,shift left=.6ex, "r"] \arrow[r,shift right=.6ex] & 2 
	\end{tikzcd}
	\]
	The genericity property does not hold for them in general. Nevertheless, we find that it does in a certain range.
	
	\begin{thm}\label{genkron}
		Assume that $r\geq 3$ and let $K_r$ be the $r$-th Kronecker quiver. If $\alpha=(a,b)$ belongs to the fundamental region of $K_r$, then the genericity property holds for $\mc{R}_{Q,\alpha}$. In particular: \[\ed_k\on{Rep}_{K_r,\alpha}\leq 1-a^2-b^2+rab+\sum_p(p^{v_p(\gcd(a,b))}-1),\]
		with equality when \Cref{collconj} holds for $d=\on{gcd}(a,b)$.
	\end{thm}
	
	\subsection*{Notational conventions}
	A base field $k$ will be fixed throughout. We will denote by $A$ an associative unital $k$-algebra. For a field extension $K/k$, we will write $A_K$ for the tensor product $A\otimes_kK$. When considering an $A_K$-module $M$, we will always assume that $M$ is a finite-dimensional $K$-vector space. For a field extension $L/K$, we will denote $M\otimes_KL$ by $M_L$.
	
	\section{Representations of quivers}\label{quiverreps}
	The purpose of this section is to briefly recall the definitions and results from the theory of quiver representations that are relevant to our discussion. 
	
	Let $Q$ be a quiver. We will write $Q_0$ for the set of its vertices and $Q_1$ for the set of its arrows between pairs of vertices. Given two vertices $i,j\in Q_0$, we will write $a:i\to j$ for an arrow having source $i$ and target $j$. When considering sums over all arrows of a quiver $Q$, the notation $\sum_{i\to j}$ will be used. We will write $\sum_{i- j}$ to indicate a sum over all edges of the underlying graph of $Q$.
	
	Let $K/k$ be a field extension. A $K$-representation $M$ of $Q$ is given by a finite-dimensional $K$-vector space $M_i$ for each vertex $i$ of $Q$, together with a linear map $\phi_a:M_i\to M_j$ for every arrow $a:i\to j$. A homomorphism of representations $f:M'\to M$ is given by $K$-linear maps $f_i:M'_i\to M_i$ such that for each arrow $a:i\to j$ one has $\phi_a\circ f_i=f_j\circ \phi'_a$. It is a basic fact that there is an equivalence of categories between $KQ$-modules and $K$-linear representations of $Q$ and that this equivalent is functorial with respect to field extensions $L/K$, see \cite[Theorem 5.4]{schiffler}. 
	
	The \emph{dimension vector} of the representation $M$ is the vector $(\dim M_i)_{i\in Q_0}$. A dimension vector $\alpha$ is said to be \emph{indivisible} if $\gcd(\alpha_i)=1$. The support of $\alpha$ is the subset $\on{supp}\alpha\c Q_0$ of vertices $i$ such that $\alpha_i\neq 0$.
	
	A quiver $Q$ is said to be of finite representation type, tame or wild if its path algebra $kQ$ is (see \cite{drozd} or \cite{scaviaalgebras} for the definitions). The connected quivers of finite representation type are classified: they are exactly those whose underlying graph is a Dynkin diagram of type $A$, $D$ or $E$. The quiver $Q$ is tame if and only if its underlying graph is an extended Dynkin diagram of type $\widetilde{A}$, $\widetilde{D}$ or $\widetilde{E}$.
	
	The Tits form of $Q$ is the bilinear form $\ang{\cdot,\cdot}:\R^{Q_0}\times\R^{Q_0}\to\R$ given by \[\ang{\alpha,\beta}:=\sum_{i\in Q_0}\alpha_i\beta_i-\sum_{i\to j}\alpha_i\beta_j.\]
	We also let $(\alpha,\beta):=\ang{\alpha,\beta}+\ang{\beta,\alpha}$.
	
	The \emph{Weyl group} of $Q$ is the subgroup $W\c\on{Aut}(\Z^{Q_0})$ generated by the \emph{simple reflections} \begin{align*}
	s_i:\Z^{Q_0}&\to\Z^{Q_0}\\
	\alpha&\mapsto \alpha-(\alpha,e_i)e_i
	\end{align*}
	where $i$ is a loop-free vertex of $Q$, and $e_i\in \Z^{Q_0}$ is the standard basis element corresponding to $i$. The \emph{fundamental region} is the set $F$ of non-zero $\alpha\in\N^{Q_0}$ with connected support and $(\alpha,e_i)\leq 0$ for all $i$. The \emph{real roots} for $Q$ are the dimension vectors that belong to an orbit of $\pm e_i$ (for $i\in Q_0$ loop-free) under the Weyl group. The \emph{imaginary roots} for $Q$ are the orbits of $\pm\alpha$ (for $\alpha\in F$) under $W$. An imaginary root $\alpha$ is called \emph{isotropic} if $\ang{\alpha,\alpha}=0$ and \emph{anisotropic} if $\ang{\alpha,\alpha}<0$. Collectively, real roots and imaginary roots are called roots.
	It can be shown that every root has either all non-negative components or all non-positive components. Hence we may speak of \emph{positive} and \emph{negative} roots.

	A dimension vector $\alpha$ is called a \emph{Schur root} if there exist a field extension $K/k$ and a $K$-representation $M$ of $Q$ of dimension vector $\alpha$ such that $\on{End}(M)=K$. Such $M$ is called a \emph{brick}.
	If $K\c L$ is a field extension and $M$ is a $K$-representation, $\on{End}_K(M)\otimes_KL=\on{End}_L(M_L)$, hence the property of being a brick is invariant under base change, and may be checked over an algebraically closed field.
	
	Given a dimension vector $\alpha$, there exists a partition $\alpha=\sum\beta_j$ such that a generic $\alpha$-dimensional representation $M$ of $Q$ is a direct sum $M=\oplus M_j$ of indecomposable representations, where $M_j$ has dimension vector $\beta_j$. For details, see \cite{kac2} and \cite{schofield91}.
	
	\section{Essential dimension of functors}\label{essdim}
	We denote by $\on{Fields}_k$ the category of field extensions of $k$. Consider a functor $F:\on{Fields}_k\to\on{Sets}$. We say that an element $\xi\in F(L)$ is \emph{defined over a field} $K\c L$, or that $K$ is a \emph{field of definition} for $\xi$, if $\xi$ belongs to the image of $F(K)\to F(L)$. The \emph{essential dimension} of $\xi$ is \[\ed_k\xi:=\min_K\trdeg_kK\] where the minimum is taken over all fields of definition $K$ of $\xi$.
	
	The \emph{essential dimension} of $F$ is defined to be
	\[\ed_kF:=\sup_{(K,\xi)}\ed_k\xi\] where the supremum is taken over all pairs $(K,\xi)$, where $K$ is a field extension of $k$, and $\xi\in F(K)$.
	
	Given a dimension vector $\alpha$, we define the functor\[\on{Rep}_{Q,\alpha}:\on{Fields}_k\to\on{Sets}\] by setting
	\[\on{Rep}_{Q,\alpha}(K):=\set{\text{Isomorphism classes of $\alpha$-dimensional $K$-representations of $Q$}}.\] If $K\c L$ is a field extension, the corresponding map $\on{Rep}_{Q,\alpha}(K)\to\on{Rep}_{Q,\alpha}(L)$ is given by tensor product.
	
	\begin{example}\label{1loop}
		Let $Q$ be the $1$-loop quiver. Then isomorphism classes of $n$-dimensional representations of $Q$ correspond to conjugacy classes of $n\times n$ matrices up to conjugation. The existence of the rational canonical form implies $\ed_k\on{Rep}_{Q,n}\leq n$. On the other hand, a matrix in rational canonical form with characteristic polynomial $t^n+a_1t^{n-1}+\dots+a_n$, with the $a_i$ algebraically independent, is defined over $k(a_1,\dots,a_n)$ but not over any proper subfield. This proves that in fact $\ed_k\on{Rep}_{Q,n}=n$. See \cite{reichsteinwhatis} for the details.
	\end{example}
	
	\begin{example}\label{deffieldrealroot}
		Let $\alpha$ be a real root for the quiver $Q$. If $K$ is an algebraically closed field, the unique indecomposable representation of dimension vector $\alpha$ is defined over the prime field of $K$. This was first proved by Kac in positive characteristic \cite[Theorem 1]{kac1} and then by Schofield in characteristic zero \cite[Theorem 8]{deffieldrealroot}. To our knowledge, it is the first result related to fields of definitions of quiver representations.
	\end{example}
	
	In \cite{bensonreichstein} and \cite{scaviaalgebras}, the following related functors are studied. Let $A$ be an associative unital $k$-algebra. For any non-negative integer $n$, we define the functor \[\on{Rep}_A[n]:\on{Fields}_k\to \on{Sets}\] by setting \[\on{Rep}_A[n](K):=\set{\text{Isomorphism classes of $n$-dimensional representations of $A_K$}}\] for every field extension $K/k$. For an inclusion $K\c L$, the corresponding map $\on{Rep}_A[n](K)\to \on{Rep}_A[n](L)$ is induced by tensor product.
	
	For a quiver $Q$, we may consider the functors $\on{Rep}_{Q,\alpha}$ for each dimension vector $\alpha$, and the functors $\on{Rep}_{kQ}[n]$ for each non-negative integer $n$. Since $K$-representations of a quiver $Q$ correspond to $K$-representations of its path algebra, functorially in $K$ there is a clear relation between the two families of functors, namely \[\ed_k\on{Rep}_{kQ}[n]=\max_{\sum\alpha_i=n}\ed_k\on{Rep}_{Q,\alpha}.\]

	\section{Essential dimension of stacks}
	
	We denote by $\on{Sch}_k$ the category of schemes over $k$. If $\mc{X}$ is an algebraic stack over $\on{Sch}_k$, we obtain a functor	\[F_{\mc{X}}:\on{Fields}_k\to\on{Sets}\] sending a field $K$ containing $k$ to the set of isomorphism classes of objects in $\mc{X}(\Spec K)$. If $\xi\in\mc{X}(K)$, we define its \emph{essential dimension} $\ed_k\xi$ to be the essential dimension of its isomorphism class in $F_{\mc{X}}$. We define the \emph{essential dimension} of $\mc{X}$ as \[\ed_k(\mc{X}):=\ed_k(F_{\mc{X}}).\] 
	
	Let $\mc{X}$ be an integral algebraic stack of finite type over a field $k$. The \emph{generic essential dimension} of $\mc{X}$ is defined as \[\on{ged}_k\mc{X}:=\sup\set{\ed_k\eta|\ \eta:\Spec K\to \mc{X} \text{ is dominant}}.\] Equivalently, it is the smallest essential dimension of a non-empty open substack of $\mc{X}$. We say that the stack $\mc{X}$ satisfies the \emph{genericity property} if \[\ed_k\mc{X}=\on{ged}_k\mc{X}.\]
	
	Let $Q$ be a quiver. It is well known that one may view $K$-representations of $Q$ as $K$-orbits of a suitable action. Let $X_{Q,\alpha}:=\prod_{i\to j}\on{Mat}_{\alpha_j\times\alpha_i,k}$ and let $G_{Q,\alpha}:=\prod_i \on{GL}_{\alpha_i,k}$ be an affine space and an algebraic group over $k$, respectively. There is an action of $G_{Q,\alpha}$ over $X_{Q,\alpha}$, given by \[(g_i)_{i\in Q_0}\cdot (P_{a})_{a:i\to j}:=(g_jP_{a}g_i^{-1})_{a:i\to j}.\] We denote by $\mc{R}_{Q,\alpha}$ the quotient stack $[X_{Q,\alpha}/G_{Q,\alpha}]$. 
	
	By \cite[Example 2.6]{genericity0}, for every field extension $K/k$, there is a natural correspondence between the orbits of this action defined over $K$, that is, $K$-points of $\mc{R}_{Q,\alpha}$, and the isomorphism classes of representations of $Q$ of dimension vector $\alpha$. Therefore \[\ed_k\on{Rep}_{Q,\alpha}=\ed_k\mc{R}_{Q,\alpha}.\]

	\begin{rmk}\label{adhoc}
		The construction of $X_{Q, \alpha}$ comes with an $\alpha$-dimensional representation $M^{\on{gen}}$ of $Q$ over the generic point $K:=k(X_{Q, \alpha})$ of $X_{Q,\alpha}$, corresponding to the natural inclusion $\Spec K\hookrightarrow X_{Q,\alpha}$. One can show that \[\on{ged}_k\mc{R}_{Q,\alpha}=\ed_kM^{\on{ged}},\] see \cite[Proposition 14.1]{berhuy2005notion}.
	\end{rmk}

	For any $k$-scheme $S$, objects of $\mc{R}_{Q,\alpha}$ over $S$ are pairs $E:=(\set{E_i}_{i\in Q_0}, \set{\phi_{a}}_{a\in Q_1})$, where $E_i$ is a locally free $\O_S$-modules of rank $\alpha_i$ for each vertex $i$ and $\phi_{a}:E_{i}\to E_{j}$ is a morphism of $\O_S$-modules for each arrow $a:i\to j$. A morphism $E'\to E$ is given by isomorphisms $E'_i\to E_i$ for each vertex $i$, satisfying the usual commutativity conditions.

	If $Q$ is a quiver of finite representation type, every root of $Q$ is a real root, so by \Cref{deffieldrealroot} every representation is defined over the prime field of $k$. In the case when $Q$ is tame, we will make use of the following observation. 
	\begin{prop}\label{tamegen}
		Let $Q$ be a tame quiver and $\delta$ be its null root. Then 
		\[\on{ged}_k\mc{R}_{Q,n\delta}=\ed_k\on{Rep}_{Q,n\delta}=n\]
	\end{prop}
	\begin{proof}
		We proved in \cite[Theorem 1.3]{scaviaalgebras} that $\ed_k\on{Rep}_{Q,n\delta}=n$ for each $n\geq 0$. Therefore, it suffices to prove that $\on{ged}_k\on{Rep}_{Q,n\delta}\geq n$. This follows from the proof of \cite[Theorem 1.3]{scaviaalgebras}, but we repeat the argument here.
		
		We may assume that $k$ is algebraically closed. There is a one-parameter family of $\delta$-dimensional indecomposable representations of $Q$. Let $Z_n\c X_{Q,n\delta}$ be the $G_{Q,\alpha}$-invariant locally closed subset parametrizing representations $\oplus_{h=1}^n M_h$, where each $M_h$ has dimension vector $\delta$. There is an obvious action of $S_n$ on $Z_n$, given by permutation of the summands. Consider $n$ copies of an infinite family of indecomposable representations of dimension vector $\delta$ parametrized by an open subset of $\A^1_k$. This gives an $S_n$-equivariant rational map \[\A^n_k\dashrightarrow Z_n\] that intersects any orbit in at most a finite number of points. By \cite[Lemma 2.3]{scaviaalgebras}, we conclude that $\ed_k\on{Rep}_{Q,n\delta}\geq n$.
	\end{proof}
	The following general lemma will be used in the proof of \Cref{genericed} and \Cref{genkron}.
	\begin{lemma}\label{allbricks}
		Let $Q$ be any quiver, $\alpha$ a Schur root for $Q$, and $M$ an $\alpha$-dimensional brick. Then \[\ed_kM\leq \on{ged}_k\mc{R}_{Q,\alpha}.\]
	\end{lemma}
	
	\begin{proof}
		Let $M$ be an $\alpha$-dimensional brick defined over $L/k$. We must show that $\ed_kM\leq \on{ged}_k\mc{R}_{Q,\alpha}$. By \Cref{candim}, this is equivalent to \[\ed_kM\leq 1-\ang{\alpha,\alpha}+\on{cd}(G_{Q,\alpha}/\mu_d),\] where $d:=\gcd(\alpha_i)$. We write $\cl{M}$ for the image of $M$ in $\cl{\mc{R}}_{Q,\alpha}$. 
		Consider a subextension $k\c K\c L$ such that $\cl{M}$ descends to $K$ and $\trdeg_kK=\ed_k\cl{M}$. We have a cartesian diagram
		\[
		\begin{tikzcd}
		\mc{G}_M \arrow[r] \arrow[d] &\mc{U}_{Q,\alpha} \arrow[d] \\
		\Spec K \arrow[r] & \cl{\mc{U}}_{Q,\alpha} 
		\end{tikzcd}
		\]\noindent
		where $\mc{G}_M$ is the residue gerbe of $M$. Since $M$ is defined over $L$, $\mc{G}_M$ is split by $L$, and the map $M:\Spec L\to \mc{U}_{Q,\alpha}$ factors through a map $M_0:\Spec L\to \mc{G}_M$. Now $M_0$ (and so $M$) descends to some intermediate subfield $K\c K_0\c L$ such that $\trdeg_kK_0\leq \ed_k(\mc{G}_M)$. By \cite[Proposition 2.3(a)]{genericity} $\ed_k\mc{G}_M=\on{cd}\mc{G}_M$. By \Cref{indgeneric} the generic gerbe has index $d$, hence it follows from \cite[Lemma 2.4(a)]{genericity} that $\on{ind}\mc{G}_M\mid d$. Therefore, by \cite[Lemma 2.2(c)]{genericity}, $\on{cd}(\mc{G}_M)\leq \on{cd}(\on{GL}_d/\mu_d)$. Consider the commutative diagram
		\[
		\begin{tikzcd}
		1\arrow[r]& \mu_d \arrow[r] \arrow[d] & \on{GL}_d  \arrow[r] \arrow[d] & \on{GL}_d/\mu_d  \arrow[r] \arrow[d] & 1\\
		1\arrow[r]& \mu_d \arrow[r] & G_{Q,\alpha} \arrow[r] & {G}_{Q,\alpha}/\mu_d \arrow[r] & 1
		\end{tikzcd}
		\]\noindent
		with exact rows. Here $\on{GL}_d$ is embedded in $G_{Q,\alpha}$ block-diagonally. The associated diagram in cohomology shows that for every field extension $k'/k$, the coboundary map $H^1(k',\on{GL}_d/\mu_d)\to H^2(k',\mu_d)=\on{Br}(k')[d]$ factors through the coboundary $H^1(k', G_{Q,\alpha}/\mu_d)\to H^2(k',\mu_d)$. Now apply \cite[Lemma 2.2(b)]{genericity} with $G=\on{GL}_d$ or $G=G_{Q,\alpha}$, and $C=\mu_d$, to obtain $\on{cd}(\on{GL}_d/\mu_d)\leq \on{cd}(G_{Q,\alpha}/\mu_d)$.
		By \cite[Lemma 2.2(c)]{genericity}, it follows that $\on{cd}(\mc{G}_M)\leq \on{cd}(G_{Q,\alpha}/\mu_d)$. On the other hand $\trdeg_kk(M)\leq \dim \cl{U}_{Q,\alpha}=1-\ang{\alpha,\alpha}$, so
		\[\ed_kM=\trdeg_kk(M)+\ed_{k(M)}M\leq 1-\ang{\alpha,\alpha}+\on{cd}(G_{Q,\alpha}/\mu_d).\] Combining this with \Cref{candim} yields $\ed_kM\leq \on{ged}\mc{R}_{Q,\alpha}$, as desired. 	
	\end{proof}
	
	\section{The Colliot-Th\'el\`ene - Karpenko - Merkurjev Conjecture}\label{section5}
	As noted in the Introduction, part of the statement of \Cref{genericed} depends on a conjecture due to Colliot-Th\'el\`ene Karpenko and Merkurjev, formulated in \cite[\S 1]{colliothelenekarpenkomerkurjev}. Following \cite{biswasdhillonhoffmann}, we rephrase this conjecture in a way that is better suited to our needs.
	
	Let $A$ be a finite-dimensional $k$-algebra. We say that a projective right $A$-module $M$ of finite dimension over $k$ has \emph{rank} $r\in\Q_{>0}$ if the direct sum $M^{\oplus n}$ is free of rank $nr$ for some $n\in\N$ with $nr\in\N$. We let $\on{Mod}_{A,r}$ be the functor of isomorphism classes of projective $A$-modules of rank $r$. 
	
	By \cite[Proposition 2.]{biswasdhillonhoffmann}, for every positive rational number $r$, $\on{Mod}_{A,r}$ is a detection functor, in the sense of \cite[\S 4a]{merkurjev2013essential}. If $A=D$ is a division algebra, and $K/k$ is a field extension, by definition $\on{Mod}_{A,r}(K)\neq \emptyset$ if and only if $X_D(K)\neq\emptyset$, where $X_D$ is the Severi-Brauer variety of $(\deg D)$-dimensional right-ideals in $D$. By \cite[\S 4a]{merkurjev2013essential}, $\ed_k\on{Mod}_{D,1/\on{deg}D}=\on{cd}(X_D)$, where $\on{cd}$ denotes the canonical dimension. We refer the reader to \cite{berhuy2005notion}, \cite{karpenko2006canonical} for an extensive treatment of the canonical dimension (denoted $\on{cd}$) of varieties and algebraic groups, and \cite[\S 2.2]{genericity} for the definition of canonical dimension of a gerbe and for a useful summary. 
	
	The following conjecture and proposition were originally stated using canonical dimension and incompressibility of $X_D$ in \cite[\S 1]{colliothelenekarpenkomerkurjev}. For our purposes, it is better to rephrase them using the functor $\on{Mod}_{D,1/\deg D}$, as is done in \cite[Conjecture 3.10]{biswasdhillonhoffmann}. 
	
	\begin{conj}\label{collconj}
		Let $d\geq 1$. If $D$ is a central division algebra of degree $d$ over $k$, then \[\ed_k(\on{Mod}_{D,1/d})=\sum_{p\mid d}(p^{v_p(d)}-1),\] the sum being over all primes $p$.
	\end{conj}
	
	\begin{prop}\label{collprop}
		Let $d\geq 1$. If $D$ is a central division algebra of degree $d$ over $k$, then \[\ed_k(\on{Mod}_{D,1/d})\leq \sum_{p\mid d}(p^{v_p(d)}-1),\] the sum being over all primes $p$. Equality holds if $d$ is a prime power or $6$.
	\end{prop}
	\begin{proof}
		The inequality is proved in \cite[\S 1]{colliothelenekarpenkomerkurjev}.
		The equality is proved in \cite[Corollary 4.4]{karpenko2013upper} when $d$ is a prime power, and in \cite[Theorem 1.3]{colliothelenekarpenkomerkurjev} when $d=6$.
	\end{proof}
	
	\section{Elementary examples}
	The following examples serve to illustrate the difference between essential dimension and generic essential dimension, in the context of quiver representations. They show that the failure of the genericity property is quite frequent. 
	
	\begin{example}
		Let $Q$ be the $2$-Kronecker quiver:
		\[
		\begin{tikzcd}[row sep=1ex]
		& 1 \arrow[r] \arrow[r,shift left=.6ex] & 2 
		\end{tikzcd}
		\]
		It is a tame quiver. The real roots of $Q$ are the dimension vectors of the form $(n,n\pm 1)$, for each $n\geq 1$. The null root of $Q$ is $\delta=(1,1)$, therefore the imaginary roots of $Q$ are of the form $n\delta=(n,n)$. Consider a real root of the form $(n,n+1)$, the other case being entirely analogous. From \Cref{deffieldrealroot}, \[\on{ged}_k\mc{R}_{Q,(n,n+1)}=0.\] On the other hand, any representation of dimension vector $(n,n+1)$ is a direct sum of indecomposable representations of dimension vector $m\delta$ or $(m,m\pm1)$. By \Cref{deffieldrealroot}, the latter are defined over the base field $k$. Using \Cref{tamegen} and \cite[Lemma 5.4]{scaviaalgebras}, we obtain 
		\[\ed_k\on{Rep}_{Q,(n,n+1)}=\ed_k\on{Rep}_{Q,n\delta}=n.\]
	\end{example}

	\begin{example}\label{starshaped}
		Let $m,n$ be non-negative integers, and consider the quiver $Q_m$ with $m+1$ vertices labeled $0,1,\dots,m$, and one arrow $a_i:i\to 0$ for every $i=1,\dots,m$. Here is a picture when $m=4$.
		\[
		\begin{tikzcd}
		& 2 \arrow[d] \\
		1 \arrow[r] & 0 & 3 \arrow[l] \\
		& 4 \arrow[u]
		\end{tikzcd}
		\]
		The quiver $Q_m$ is of finite representation type when $m\leq 3$, tame when $m=4$, and wild for $m\geq 5$. As dimension vector, choose $\alpha_{m,n}:=(n+1,1,\dots,1)$. An $\alpha_{m,n}$-dimensional representation of $Q_m$ over $K$ is given by at most $m$ lines in $K^{n+1}$, up to linear automorphisms of $K^{n+1}$. It is basically the datum of at most $m$ points in $\P^{n}_K$ up to projective equivalence. More precisely, consider the functor $\on{Rep}_{Q_m,\alpha_{m,n}}$ and \begin{align*}
		L_{m,n}:\on{Fields}_k&\to \on{Sets}\\
		K&\mapsto\set{\text{$\on{PGL}_{n+1}$-orbits in $(\P^n\cup\set{0})^m(K)$}}
		\end{align*}
		where $\on{PGL}_{n+1}$ acts diagonally on $(\P^n)^r$ for every $0\leq r\leq m$, and fixes $0$. 
		
		There is a morphism of functors $\Phi:\on{Rep}_{Q_m,\alpha_{m,n}}\to L_{m,n}$ constructed as follows. If $(M,\phi)$ is a $K$-representation, fix an isomorphism $\P(M_0)\cong \P^n_K$. Then $(M,\phi)$ gets sent to the $K$-point of $(\P^n)^m$ whose $r$-th component is $\on{Im}\phi_{\alpha_i}$ when it is not zero, and the point $0$ otherwise. Of course, the orbit associated to $(M,\phi)$ in this way does not depend on the choice of an isomorphism $M_0\cong K^{n+1}$.
		
		We want to show that $\Phi$ is an isomorphism. It is immediate to check that if two $K$-representations map to the same orbit, then they are isomorphic, so $\Phi$ is injective. Given a $K$-orbit $\mc{O}$ of $(\P^n\cup\set{0})^m$, choose a $K$-point $(L_1,\dots, L_m)\in \O$. Set $M_0:= K^n, M_i:=K$ for $i\geq 1$ and let $\phi_{\alpha_i}$ be the zero map if $L_i=0$, and send $1$ to any non-zero vector lying on the line $L_i$ otherwise. This defines a representation $(M,\phi)$ such that $\Phi(M,\phi)=\O$ so $\Phi$ is surjective. Hence $\Phi$ is an isomorphism. In particular, $\ed_k\on{Rep}_{Q_m,\alpha_{m,n}}=\ed_kL_{m,n}$.

		We start by computing the generic essential dimension of $\on{Rep}_{Q_m,\alpha_{m,n}}$. If a map $\Spec K\to \mc{R}_{Q_m,\alpha_{m,n}}$ is dominant, the corresponding orbit in $(\P^n\cup\set{0})^m(K)$ consists of $m$-uples of points in $(\P^n)^m$ in general position. If $m\leq n+2$ then $\on{PGL}_{n+1}$ acts transitively on $m$-uples of points in general position. If $m>n+2$ and the points are in general position, we may assume after acting with $\on{PGL}_{n+1}$ that $n+2$ of them will be of the form \begin{equation}\label{genpos}(1:0:\dots:0),(0:1:0:\dots:0),\dots,(0:\dots:0:1),(1:\dots:1).\end{equation} The $\on{PGL}_{n+1}$-orbit of this $m$-tuple is then completely determined by the remaining $m-n-2$ points. Since any one of them is determined by $n+1$ coordinates up to simultaneous rescaling, each of the $m-n-2$ points contributes at most $n$ to the essential dimension. Moreover, consider the configuration of $m$ points, where the first $n+2$ are as in (\ref{genpos}), and the remaining $m-n-2$ are of the form \[(1:a_{i1}:\dots:a_{in}), \quad i=1,\dots,m-n-2,\] where the $a_{ij}$ are independent variables over $k$. This configuration has a minimal field of definition $K:=k(a_{ij})_{i,j}$, so that $\trdeg_kK=n(m-n-2)$. Moreover, the corresponding map $\Spec K\to \mc{R}_{Q_m,\alpha_{m,n}}$ is dominant. We obtain:
		\[
		\on{ged}_k\mc{R}_{Q_m,\alpha_{m,n}}=\begin{cases}
		0 \text{ if $m\leq n+2$,}\\
		n(m-n-2) \text{ if $m>n+2$.}
		\end{cases}
		\]
		We now determine the essential dimension of $R_{Q_m,\alpha_{m,n}}$. In order to compute it, we may clearly restrict ourselves to representations $(M,\phi)$ such that $\phi_{\alpha_i}\neq 0$ for every $i$, that is, $\on{PGL}_{n+1}$-orbits in $(\P^n)^m$. Consider a configuration of points spanning a subspace $H$ of $\P^n$ of dimension $r\leq \min(n,m-1)$. After a translation by an element of $\on{PGL}_{n+1}$, we may assume that $H$ is given by the vanishing of the last $n-r$ coordinates. If $m=r+1$, $\on{PGL}_{n+1}$ acts transitively on $m$-uples of points of $H$. If $m\geq r+2$, the action of $\on{PGL}_{n+1}$ may be used to put $r+2$ points in the form \[(1:0:\dots:0),\dots,(0:\dots:0:1:0:\dots:0), (1:\dots:1:0:\dots:0).\] 
		The remaining $m-r-2$ points are now fixed, and are determined by $r+1$ coordinates up to scaling. Using the inequality $ab\leq \frac{1}{4}(a+b)^2$, it is easy to see that the essential dimension of $R_{Q_m,\alpha_{m,n}}$ is at most:
		\begin{equation}\label{esempioed}
		\max_{1\leq r\leq\min(n,m-1)}r(m-r-2)=\begin{cases}
		\frac{1}{4}(m-2)^2 \text{ if $m\leq 2n$ is even,}\\
		\frac{1}{4}(m-1)(m-3) \text{ if $m\leq 2n$ is odd,}\\
		n(m-n-2) \text{ if $m>2n$.}
		\end{cases}
		\end{equation}
		Moreover, one can construct examples showing that equality actually holds, in a way which is totally analogous to what we did for $\on{ged}_k\mc{R}_{Q_m,\alpha_{m,n}}$, so $\ed_kR_{Q_m,\alpha_{m,n}}$ is given by (\ref{esempioed}).
		
		This gives a very explicit class of examples for which the genericity property does not hold. The simplest among these examples is when $m=4$ and $n=2$. In this case $Q=\widetilde{D}_4$ is tame, and $\alpha_{4,2}=(3,1,1,1,1)$. Since $\on{PGL}_3$ acts transitively on $4$-uples of points in $\P^2$ in general position, the generic essential dimension is zero. On the other hand, if the $4$ points lie on a common line, the essential dimension may be $1$.
	\end{example}

	\section{Proof of \Cref{genericed}}
	
	Let $\mc{X}$ be an irreducible algebraic stack. Then $\mc{X}$ admits a \emph{generic gerbe}, defined as the residual gerbe at any dominant point $\Spec K\to\mc{X}$ (see \cite[Chapter 11]{laumonmoretbailly}). If $\alpha$ is a Schur root for the quiver $Q$, the generic $\alpha$-dimensional representation is a $\Gm$-gerbe, and so gives rise to a Brauer class in $\on{Br}(k(\mc{G}))$. In order to understand the essential dimension of the generic gerbe of $\mc{R}_{Q,\alpha}$, the first step is to compute its index.
	
	\begin{lemma}\label{residuegerbe}
		Let $\mc{G}$ be the residue gerbe of a brick of $\mc{R}_{Q,\alpha}$. Then $\on{ind}\mc{G}\mid\gcd(\alpha_i)$.
	\end{lemma}  
	Recall that a brick is a $K$-representation $M$ of $Q$ such that $\on{End}_KM=K$.
	\begin{proof}
		Since $\mc{G}$ parametrizes bricks, it is a $\Gm$-gerbe, so its index is well-defined. By \cite[Lemma 3.10]{hoffmann} we know that $\ind\mc{G}$ is the greatest common divisor of the ranks of all the twisted sheaves (i.e., vector bundles of rank $1$) on some open substack of $\mc{R}_{Q,\alpha}$. 
		
		To prove that $\ind\mc{G}\mid \gcd(\alpha_i)$, it is therefore sufficient to exhibit for every $i\in Q_0$ a twisted sheaf on $\mc{R}_{Q,\alpha}$ of rank $\alpha_i$. Recall that a vector bundle of rank $r$ on $\mc{R}_{Q,\alpha}$ is a $1$-morphism $\mc{V}:\mc{R}_{Q,\alpha}\to\mc{V}ect_{r}$. If $S$ is a scheme over $k$, an object of $\mc{R}_{Q,\alpha}(S)$ is a pair $E:=(\set{E_i}_{i\in Q_0}, \set{\phi_{a}}_{a\in Q_1})$, where $E_i$ is a vector bundle over $S$ of rank $\alpha_i$ for each vertex $i$ and $\phi_{a}:E_{i}\to E_{j}$ is a morphism $\O_S$-modules for each arrow $a:i\to j$. Fix a vertex $i_0\in Q_0$, and set $\mc{V}(E):=E_{i_0}$. Now let $E\in\mc{R}(S)$ and $E':=(\set{E'_i}_{i\in Q_0}, \set{\phi'_{a}}_{a\in Q_1})\in\mc{R}(S')$, where $S'$ is also a scheme over $k$ and let $f:=(f_i:E'_i\to E_i)_{i\in Q_0}$ be a morphism from $E'$ to $E$ in $\mc{R}$, set $\mc{V}(f):=f_{i_0}$. By definition, $\mc{V}$ is a vector bundle of weight $1$ and rank $\alpha_{i_0}$.
	\end{proof}

	\begin{prop}\label{indgeneric}
		Let $\alpha$ be a Schur root. The index of the generic gerbe of $\mc{R}_{Q,\alpha}$ is equal to $\gcd(\alpha_i)$.	
	\end{prop}
	
	\begin{proof}
		Let us call $\mc{G}$ the generic gerbe of $\mc{R}_{Q,\alpha}$. By \Cref{residuegerbe}, we have $\on{ind}\mc{G}\mid \gcd(\alpha_i)$, so it suffices to show that $\gcd(\alpha_i)|\ind\mc{G}$. The proof will follow from the next two lemmas.
		
		\begin{lemma}
			Suppose that there is a line bundle $\mc{L}$ of weight $w\in\Z$ on an open substack $\mc{U}$ of $\mc{R}_{Q,\alpha}$. Then we may extend $\mc{L}$ to a line bundle $\mc{L}'$ on $\mc{R}_{Q,\alpha}$ of the same weight.
		\end{lemma}
		
		\begin{proof}
			We make use of the following standard result.
			
			\begin{fact}
				Let $\mc{X}$ be a noetherian algebraic stack over $k$ and $\mc{U}$ an open substack of $\mc{X}$. Denote by $j:\mc{U}\to\mc{X}$ the inclusion $1$-morphism. Let $\mc{M}$ be a quasi-coherent $\O_{\mc{X}}$-module and $\mc{N}$ a coherent $\O_{\mc{U}}$-submodule of $j^*\mc{M}$. Then there exists a coherent $\O_{\mc{X}}$-submodule $\mc{N}'$ of $\mc{M}$ such that $j^*\mc{N}'=\mc{N}$.
			\end{fact}
			
			\begin{proof}
				See \cite[Corollaire 15.5]{laumonmoretbailly}.
			\end{proof}
			
			In our case we take $\mc{X}=\mc{R}_{Q,\alpha}$, $\mc{N}=\mc{L}$ and $\mc{M}=j_*\mc{L}$. Since $\mc{R}_{Q,\alpha}$ is noetherian, $\mc{M}$ is quasi-coherent. The lemma gives us a coherent subsheaf $\mc{F}\c j_*\mc{L}$. Then the double dual $\mc{L}':=\mc{F}^{**}$ is a reflexive coherent sheaf of rank one on a smooth stack, hence a line bundle (this follows immediately by \cite[VII 4.2]{bourbakialgebra}). The weight of $\mc{L}'$ is $w$ because this may be checked on $\mc{U}$, where $\mc{L}'$ restricts to $\mc{L}$. 
		\end{proof}

		\begin{lemma}\label{index2}
			Suppose that there is a line bundle $\mc{L}$ of weight $w\in\Z$ on an open substack $\mc{U}$ of $\mc{R}_{Q,\alpha}$. Then $\gcd(\alpha_i)\mid w$.
		\end{lemma}
		
		\begin{proof}
			By the previous lemma we may assume that $\mc{L}$ is defined on $\mc{R}_{Q,\alpha}$. 
			Denote by $S_{\alpha}\in \mc{R}_{Q,\alpha}(k)$ the trivial representation of $Q$ of dimension vector $\alpha$ over $k$, for which the linear maps are all zero. Then the central $\Gm\c \on{GL}_{\alpha}:=\prod_i \on{GL}_{\alpha_i}=\on{Aut}(S_{\alpha})$ acts with weight $w$ on the fiber of $\mc{L}'$ over $S_{\alpha}$. Since any one-dimensional representation of $\on{GL}_\alpha$ is of the form \[(A_1,\dots,A_r)\mapsto \det(A_1)^{m_1}\cdot\ldots\cdot\det(A_r)^{m_r}\]
			we get $w=m_1\alpha_1+\dots+m_r\alpha_r$, by restricting the above formula to $r$-uples of diagonal matrices. Hence $w$ is a multiple of $\gcd(\alpha_i)$.
		\end{proof}
		
		We are now ready to complete the proof of \Cref{indgeneric}. Let $\mc{V}$ be a vector bundle of rank $n$ and weight $w$ on some open substack $\mc{U}$. Define $\mc{L}:=\on{det}(\mc{M})$, then $\mc{L}$ is a line bundle of weight $nw$. In particular, if $\mc{V}$ has weight $1$, $\mc{L}$ has weight $n$, so by \Cref{index2} $\gcd(\alpha_i)\mid n$. Hence $\gcd(\alpha_i)\mid\ind\mc{G}$, as wanted.
	\end{proof}
	
	Let $\mc{G}$ be the residue gerbe of a brick in $\mc{R}_{Q,\alpha}$, for some Schur root $\alpha$. Since $\mc{G}$ parametrizes bricks, it is a $\Gm$-gerbe, and so admits a Brauer class in $\on{Br}(k(\mc{G}))$. On the other hand, by the Nullstellensatz, there exists a field extension $l/k(\mc{G})$ of finite degree $d$ such that $\mc{G}(l)$ is non-empty. If $V\in\mc{G}(l)$,  \[R:=\on{End}_{k(\mc{G})}(V)\] is a central simple algebra over $k(\mc{G})$ split by $l$. It is not hard to check that this class is independent of the chosen field extension $l/k(\mc{G})$.
	
	\begin{lemma}\label{classescoincide}
		The Brauer classes of $\mc{G}$ and $R$ in $\on{Br}(k(\mc{G}))$ coincide.
	\end{lemma}
	\begin{proof}
		We briefly recall the construction of the Brauer class of $\mc{G}$, as given in \cite[Lemma 3.10]{hoffmann}. One starts by choosing a field extension $l/k(\mc{G})$ of finite degree $d$ such that $\mc{G}(l)$ is non-empty. This means that $\mc{G}_{l}\cong B\Gm$, so it admits a line bundle $\mc{L}_1$ of weight $1$, corresponding to the tautological $1$-dimensional representation of $\Gm$. If $\pi:\mc{G}_{l}\to\mc{G}$ denotes the natural projection, $\mc{V}:=\pi_*\mc{L}_1$ is a vector bundle of rank $d$ and weight $1$ on $\mc{G}$. The algebra bundle $\mc{E}nd(\mc{V})$ on $\mc{G}$ has weight $0$, and so descends to a central simple algebra $A$ split by $l$. By definition, the Brauer class of $\mc{G}$ is that of $A$. One then checks that this definition does not depend on the choice of the extension $l/k(\mc{G})$. 	
		
		Moreover, there is a chain of isomorphisms of $k(\mc{G})$-vector spaces.
		\begin{align*}
		R=\on{Hom}_{k(\mc{G})}({V},{V})&\cong \on{Hom}_l({V}\otimes_{k(\mc{G})}l,V)\\
		&\cong\on{Hom}_{l}(V^d,V)\\
		&\cong\on{Hom}_{l}(V,V)^d\\
		&\cong\on{Hom}_{l}(\mc{L}_1(V),\mc{L}_1(V))^d\\
		&\cong\on{Hom}_{l}(\mc{L}_1(V)^d,\mc{L}_1(V))\\
		&\cong\on{Hom}_{l}(\pi^*\mc{V}(V),\mc{L}_1(V))\\
		&\cong\on{Hom}_{k(\mc{G})}(\mc{V}(V),\mc{V}(V))=A.
		\end{align*}	
		The map $\on{Hom}_{l}(V,V)\to\on{Hom}_{l}(\mc{L}_1(V),\mc{L}_1(V))$ is the one induced by the functor $\mc{L}$. The map $\on{Hom}_{k(\mc{G})}({V},{V})\to\on{Hom}_{k(\mc{G})}(\mc{V}(V),\mc{V}(V))$ is exactly the map given by the functor $\mc{V}$, hence also respects composition. Thus it is an isomorphism of $k(\mc{G})$-algebras.
	\end{proof}
	
	\begin{proof}[Proof of \Cref{genericed}]
		Let $\on{Spec}K\to \mc{R}_{Q,\alpha}$ be a dominant map, corresponding to an $\alpha$-dimensional $K$-representation $M$. Since $\alpha$ is a Schur root, $M$ is a brick. Let $M$ be a brick of dimension $\alpha$. We have \[\trdeg_kk(M)\leq \dim\on{Aut}(M)+\dim\mc{R}_{Q,\alpha}=1-\ang{\alpha,\alpha}.\] Let $\mc{G}$ be the residue gerbe of $M$. It is a $\Gm$-gerbe with residue field $k(M)$. From \cite[Theorem 6.2]{scaviaalgebras} we see that \[\ed_{k(M)}\mc{G}=\ed_{k(M)}(\on{Mod}_{R,1/\on{deg}R}),\] for some central simple algebra $R$ over $k(M)$ split by $l$. By \Cref{classescoincide} and \Cref{residuegerbe}, the index of $R$ divides $\gcd(\alpha_i)$. Inequality (\ref{coll1}) now follows from \cite[Corollary 3.8]{biswasdhillonhoffmann}. 
		Furthermore, by \Cref{indgeneric} the index of the generic gerbe is $\gcd(\alpha_i)$, so equality in (\ref{coll1}) follows from \Cref{collconj} and \Cref{collprop}, for $d=\gcd(\alpha_i)$.
	\end{proof}
	The following is a corollary of the proof of \Cref{genericed}.

	\begin{cor}\label{arbitraryroot}
		Let $\alpha$ be a root of $Q$. If the canonical decomposition of $\alpha$ consists only of real roots, then \[\on{ged}_k\mc{R}_{Q,\alpha}=0.\] Otherwise, let $\beta$ be the unique imaginary Schur root appearing in the canonical decomposition of $\alpha$; see \cite[Theorem 4.4]{schofield91}. If $\beta$ is isotropic of multiplicity $m\geq 1$, then \[\on{ged}_k\mc{R}_{Q,\alpha}=m.\]
		If $\beta$ is anisotropic, then  \begin{equation}\label{coll2}\on{ged}_k\mc{R}_{Q,\alpha}\leq 1-\ang{\beta,\beta}+\sum_{p}(p^{v_p(\gcd(\beta_i))}-1).\end{equation}
		One has equality if \Cref{collconj} holds for $d=\gcd(\alpha_i)$.
	\end{cor}	
	
	\begin{proof}
		Our argument will make use of the reflection functors. We refer the reader to \cite[Section 3.2]{kirillov} for background material on reflection functors. We note that reflection functors may be defined over any field, and their formation commutes with extension of scalars. It is an immediate consequence of \cite[Theorem 3.11]{kirillov} that if $\sigma_i$ is a reflection at an admissible vertex $i$ (a source or a sink), and $\alpha$ is a Schur root, then \[\on{ged}_k\mc{R}_{Q,\alpha}=\on{ged}_k\mc{R}_{Q',\sigma_i(\alpha)}\] where $Q'$ is obtained from $Q$ by reversing all the arrows at $i$.
		
		Let now $\alpha$ be a root. By \cite[Theorem 4.4]{schofield91} the canonical decomposition of $\alpha$ contains at most one imaginary root. If all roots are real, by \Cref{deffieldrealroot} the generic representation is a direct sum of indecomposable representations, all of which are defined over the prime field of $k$ by \Cref{deffieldrealroot}. Hence $\on{ged}_k\on{Rep}_{Q,\alpha}=0$ in this case. Assume now that there exists an imaginary root $\beta$ in the canonical decomposition of $\alpha$, and let $M$ be a generic $\alpha$-dimensional representation. 
		
		Using a suitable sequence of reflection functors we may assume that $\beta$ is in the fundamental region of $Q$. We remark that although the reflection functors change orientation of the arrows, the fundamental region does not change. By \cite[Proposition 4.14]{lebruyn}, $\beta$ is either an anisotropic root Schur root, or is a multiple of the null root of some tame subquiver of $Q$. In the first case, one may apply \cite[Lemma 5.4]{scaviaalgebras} and the first part of the theorem to conclude. In the second case, the result follows from \cite[Lemma 5.4]{scaviaalgebras} and \Cref{tamegen}.
	\end{proof}
	
	\begin{rmk}\label{edp}
		If we consider generic essential $p$-dimension (see \cite[\S 1.1]{merkurjev2013essential}), the inequalities (\ref{coll1}) and (\ref{coll2}) of \Cref{genericed} become unconditional equalities:
		\[\on{ged}_{k,p}\mc{R}_{Q,\alpha}=1-\ang{\alpha,\alpha}+\max_{p}(p^{v_p(\gcd(\alpha_i))}-1)\]
		and \[\on{ged}_{k,p}\mc{R}_{Q,\alpha}=1-\ang{\beta,\beta}+\max_{p}(p^{v_p(\gcd(\beta_i))}-1).\]
		In the notation of (\ref{coll1}) and (\ref{coll2}), this gives unconditional lower bounds
		\[\on{ged}_k\mc{R}_{Q,\alpha}\geq 1-\ang{\alpha,\alpha}+\max_p(p^{v_p(\gcd(\alpha_i))}-1)\]
		and \[\on{ged}_{k}\mc{R}_{Q,\alpha}\geq 1-\ang{\beta,\beta}+\max_{p}(p^{v_p(\gcd(\beta_i))}-1).\]
	\end{rmk}
	
	We now give an unconditional formula for $\on{ged}_k\mc{R}_{Q,\alpha}$, involving canonical dimension; see \Cref{section5} for references on this notion.
	
	\begin{prop}\label{candim}
		Let $\alpha$ be a Schur root for the quiver $Q$, and set $d:=\gcd(\alpha_i)$. Then \begin{equation}\on{ged}_k\mc{R}_{Q,\alpha}=1-\ang{\alpha,\alpha}+\on{cd}(G_{Q,\alpha}/\mu_d).\end{equation}
	\end{prop}
	Recall that $G_{Q,\alpha}:=\prod_i \on{GL}_{\alpha_i,k}$, the product being over all $i\in Q_0$. Here $\mu_d$ is embedded in $G_{Q,\alpha}$ as the subgroup $\set{(\zeta \on{Id}_{\alpha_i})_{i\in Q_0}:\zeta^d=1}$. We will use \Cref{candim} in the proof of \Cref{genericityproperty}.
	\begin{proof}
		This argument was inspired, in part, by the proof of \cite[Proposition 7.1]{genericity}. Let $\cl{G}_{Q,\alpha}:=G_{Q,\alpha}/H$, where $H\cong \G_m$ is the diagonal copy of $\G_m$ inside ${G}_{Q,\alpha}$, and set $\cl{\mc{R}}_{Q,\alpha}:=[X_{Q,\alpha}/\cl{G}_{Q,\alpha}]$. Let $U_{Q,\alpha}$ be the open subscheme of $X_{Q,\alpha}$ parametrizing bricks, and define the open substacks $\mc{U}_{Q,\alpha}:=[U_{Q,\alpha}/G_{Q,\alpha}]$ and $\cl{\mc{U}}_{Q,\alpha}:=[U_{Q,\alpha}/\cl{G}_{Q,\alpha}]$ of $\mc{R}_{Q,\alpha}$ and $\cl{\mc{R}}_{Q,\alpha}$, respectively. We have a cartesian diagram
		\[
		\begin{tikzcd}
		\mc{U}_{Q,\alpha} \arrow[r] \arrow[d] & \mc{R}_{Q,\alpha}\arrow[d,"\pi"] \\
		\cl{\mc{U}}_{Q,\alpha} \arrow[r]
		& \cl{\mc{R}}_{Q,\alpha} 
		\end{tikzcd}
		\]\noindent
		where the horizontal maps are open embeddings, and the vertical maps are $\G_m$-gerbes. 
		Since $\alpha$ is a Schur root, $U_{Q,\alpha}$, $\mc{U}_{Q,\alpha}$ and $\cl{\mc{U}}_{Q,\alpha}$ are non-empty. Moreover, $G_{Q,\alpha}$ acts freely on $U_{Q,\alpha}$, so $\cl{\mc{U}}_{Q,\alpha}$ is an integral algebraic space of finite type. It has dimension $1-\ang{\alpha,\alpha}$; see \cite[Proposition 4.4]{king}. We set $d:=\gcd(\alpha_i)$.
		
		Let $\mc{G}$ be the generic gerbe of $\mc{R}_{Q,\alpha}$, i.e. the generic fiber of $\pi$. Its residue field is $k(\mc{G}):=k(\cl{\mc{U}}_{Q,\alpha})$. Then $\on{ged}_k\mc{R}_{Q,\alpha}=1-\ang{\alpha,\alpha}+\ed_{k(\mc{G})}\mc{G}$. If $\gamma$ denotes the class of $\mc{G}$ in $H^2(k(\mc{G}),\mc{G})$, then by \cite[Proposition 2.3(a)]{genericity} $\ed_{k(\mc{G})}\mc{G}=\on{cd}\gamma$. 
		
		The action of $\cl{G}_{Q,\alpha}$ on $X_{Q,\alpha}$ is linear and generically free, hence it gives rise to a versal $\cl{G}_{Q,\alpha}$-torsor $t\in H^1(k(\mc{G}),\cl{G}_{Q,\alpha})$, and $\gamma$ is the image of $t$ under the boundary map $H^1(k(\mc{G}), \cl{G}_{Q,\alpha})\to H^2(k(\mc{G}), \G_m)$ associated with the exact sequence \[1\to \G_m \to G_{Q,\alpha}\to \cl{G}_{Q,\alpha}\to 1.\] Since $t$ is versal, $\on{cd}t=\on{cd}{\cl{G}_{Q,\alpha}}$; see \cite[\S 2.2]{genericity}. On the other hand, by \cite[Lemma 2.2(b)]{genericity} $\on{cd}t=\on{cd}{\gamma}$. By \cite[Corollary A.2]{cernele2015essential} $H^1(-, \cl{G}_{Q,\alpha})=H^1(-,G_{Q,\alpha}/\mu_d)$, hence $\on{cd}(\cl{G}_{Q,\alpha})=\on{cd}(G_{Q,\alpha}/\mu_d)$. Combining these equalities we obtain \[\ed_{k(\mc{G})}(\mc{G})=\on{cd}\gamma=\on{cd}t=\on{cd}{\cl{G}_{Q,\alpha}}=\on{cd}(G_{Q,\alpha}/\mu_d).\]  This proves \Cref{candim}.
	\end{proof}
	
	\section{Fields of definition}

	If $M$ is a representation of $Q$, we denote by $k(M)$ its residue field, i.e., the residue field of its residue gerbe (see \cite[Chapter 11]{laumonmoretbailly}). Since $k(M)$ contains any field of definition for $M$, we have \[\ed_kM=\ed_{k(M)}M+\trdeg_kk(M).\] In this section, we address the first term of this sum, by presenting a strengthening of \cite[Lemma 6.6]{scaviaalgebras} for quiver algebras.
	
	\begin{lemma}\label{finitefieldext}
		Let $M$ be a representation of $Q$, and let $\mc{G}$ be its residue gerbe in $\mc{R}_{Q,\alpha}$, with residue field $K:=k(\mc{G})$.
		There exists a separable finite field extension $L$ of $K$ such that $\mc{G}(L)\neq\emptyset$.
	\end{lemma}
	
	\begin{proof}
		Let $\alpha$ be the dimension vector of $M$. Since $\mc{R}_{Q,\alpha}$ is of finite type over $k$, by \cite[Th?or?me 11.3]{laumonmoretbailly} the gerbe $\mc{G}$ is of finite type over $K$. We may find a smooth cover $U\to\mc{G}$ that is of finite type over $K$. Let $\Spec l\to U$ be a closed point. Then, by the Nullstellensatz, $l$ is a finite extension of $K$. The composition $\Spec l\to U\to\mc{G}$ gives an $l$-point for $\mc{G}$, corresponding to an object $\xi\in\mc{G}(l)$. This is equivalent to $\mc{G}_l\cong {B}\on{Aut}(\xi)$. Since $\on{Aut}(\xi)$ is an open subscheme of a vector space, it is smooth, hence $\mc{G}$ is smooth. Hence $U$ is also smooth over $k$. It follows that there is a closed point of $U$ whose residue field $L$ is separable over $K$, and this gives an $L$-point of $\mc{G}$.
	\end{proof}	
	
	\begin{prop}\label{estimate}
		Let $Q$ be a quiver, and let $M$ be an indecomposable $\alpha$-dimensional $K$-representation of $Q$, for some field $K$ containing $k$. Then \[\ed_{k(M)}M\leq \min_{i\in \on{supp}\alpha}\alpha_i-1.\]
	\end{prop}
	
	\begin{proof}
		Let $\mc{G}$ be the residue gerbe of the point $\Spec K\to \mc{R}_{Q,\alpha}$ given by $M$. By \Cref{finitefieldext} there exist a separable finite field extension $l$ of the residue field $k(\mc{G})=k(M)$ and an $l$-representation $N$ of $Q$ such that $N_{L}\cong M_{L}$ for any field $L$ containing both $K$ and $l$. We may assume that $l/k(M)$ is Galois with Galois group $G$, and we set $d:=[l/k(M)]$. We denote by $\cl{N}$ the $k(M)$-representation of $Q$ obtained from $N$ by restriction of scalars. Let \[\cl{N}=\oplus_{h=1}^s N_h^{\oplus r_h}\] be the decomposition of $\cl{N}$ in indecomposable $k(M)$-representations of $Q$, where $N_h\cong N_{h'}$ if and only if $h=h'$.
		Let $L=lK$ be a compositum of $l$ and $K$. Since $M$ is defined over $K$ and is indecomposable, the Galois group of $L/K$ acts transitively on isomorphism classes of indecomposable summands of $M_L$, so the $N_h$ all have the same dimension vector $\beta$. In particular, $\on{supp}\beta=\on{supp}\alpha$. By definition
		\[d\alpha=\dim_{k(M)}\cl{N}=\sum_h r_h\dim_{k(M)}N_h=(\sum_h r_hn_h)\beta.\]
		On the other hand, we may write \[\cl{N}\otimes_{k(M)}l=\oplus_{\sigma\in G}N^{\sigma}.\] Since all indecomposable summands of $N_L$ have dimension vector $\beta$, the same is true for those of $N^\sigma_L$. Since every $N_h$ is a summand of $\cl{N}$, we have just shown that each $N_h$ has dimension vector multiple of $\beta$. For every $h$, we write $\dim_{k(M)}N_h=n_h\beta$, where $n_h\geq 1$. 
		
		Consider \[A:=\on{End}_{k(M)}(N)/j(\on{End}_{k(M)}(N))\] and \[A_h:=\on{End}_{k(M)}(N_h)/j(\on{End}_{k(M)}(N_h))\] for $h=1,\dots,s$. We may write $A_h=M_{r_h}(D_h)$ for some division algebra $D_h$. Fitting's lemma and \cite[Corollary 3.7]{biswasdhillonhoffmann} imply
		\[A=\prod_{h=1}^sA_h.\]
		Let $i\in\on{supp}\alpha$. By \cite[Lemma 6.5]{scaviaalgebras} $\dim_{k(M)}D_h\leq \dim(N_i)_h$ for every $h$. By \cite[Corollary 3.7]{biswasdhillonhoffmann} \[\ed_{k(M)}(\on{Mod}_{A_h,1/d})<\frac{r_h}{d}\dim_{k(M)}(N_h)_i.\] 
		Using \cite[Proposition 3.3 and Proposition 3.2]{biswasdhillonhoffmann}, we get
		\[\ed_{k(M)}(\on{Mod}_{A,1/d})\leq \sum_h\ed_{k(M)}(\on{Mod}_{A_h,1/d}) <\frac{1}{d}\sum_hr_h\dim_{k(M)}(N_h)_i=\alpha_i\] for each vertex $i\in \on{supp}\alpha$. The claimed inequality now follows from an application of \cite[Theorem 6.2]{scaviaalgebras}.
	\end{proof}

	\section{Beginning of the proof of \Cref{genericityproperty}}
	\label{section9}
	In this section we show that if the quiver $Q$ is of finite representation type or admits at least one loop at every vertex, then $\mc{R}_{Q,\alpha}$ satisfies the genericity property for every dimension vector $\alpha$. This will establish one direction of \Cref{genericityproperty}; we will prove the other direction in the next section. 
	
	By \cite[Remark 9.1]{scaviaalgebras}, if $Q$ is of finite representation type, \begin{equation}\label{frt}\ed_k\on{Rep}_{Q,\alpha}=0\end{equation} for every dimension vector $\alpha$, so the genericity property holds in this case.

	Consider now the case where $Q$ has at least one loop at every vertex. We start by reducing the problem to the following assertion. Recall that a dimension vector $\alpha$ for $Q$ is called sincere if $\alpha_i\neq 0$ for every $i\in Q_0$.
	
	\begin{claim}\label{claim}
		Let $Q$ be a quiver having at least one loop at every vertex. Assume that $Q$ is not the $1$-loop quiver. Then for every sincere dimension vector $\alpha$, and for every $\alpha$-dimensional representation $M$ of $Q$ that is not a brick, we have \[\ed_kM\leq -\ang{\alpha,\alpha}.\]
	\end{claim}
	
	\begin{lemma}\label{reduction1}
		Assume that \Cref{claim} holds for every quiver having at least one loop at every vertex. Let $Q$ be a quiver with at least one loop at every vertex, and let $\alpha$ be a dimension vector for $Q$. Then $\mc{R}_{Q,\alpha}$ satisfies the genericity property.
	\end{lemma}
	
	\begin{proof}
		Since $Q$ has at least one loop at every vertex, every dimension vector $\alpha$ belongs to the fundamental region, hence by \cite[Proposition 4.14]{lebruyn} either it has tame support or is an imaginary anisotropic Schur root. On the other hand, the only tame quiver with at least one loop at every vertex is the $1$-loop quiver, so if $\alpha$ has tame support then $\alpha=me_i$ for some $m\geq 1$ and some vertex $i$. For such $\alpha$ the genericity property is easily seen to be true (see \Cref{1loop}). Assume now that $\alpha$ is an imaginary anisotropic Schur root. The subquiver $Q'$ of $Q$ defined by $Q'_0=\on{supp}\alpha$ and $Q'_1$ the set of all arrows in $Q_1$ between vertices in $\on{supp}\alpha$ also has one loop at each vertex, thus we are reduced to the case when $\alpha$ is sincere.
		
		When $\alpha$ is sincere, by \Cref{claim} $\ed_kM\leq -\ang{\alpha,\alpha}$ for every representation $M$ that is not a brick. By \Cref{edp} (or more directly by \cite[Lemma 10.1]{scaviaalgebras}), $\on{ged}_k\mc{R}_{Q,\alpha}\geq 1-\ang{\alpha,\alpha}$, so the maximum must be attained among bricks. The conclusion now follows from \Cref{allbricks}.	 
	\end{proof}

	\begin{prop}\label{reduction2}
		Let $Q$ be a quiver with at least one loop at every vertex and that is not the $1$-loop quiver. Then \Cref{claim} holds for $Q$.
	\end{prop}
	The combination of \Cref{reduction1} and \Cref{reduction2} proves the first implication of \Cref{genericityproperty}.
	
	\begin{proof}
		Let $\alpha$ be a sincere dimension vector for $Q$. We must show that for every field extension $K/k$ and every $K$-representation $M$ of $Q$ \[\ed_kM\leq -\ang{\alpha,\alpha}.\]	
		
		For each vertex $i$ of $Q$, let $l_i$ be the number of loops at $i$. Since $Q$ has at least one loop at every vertex, we have $l_i\geq 1$ for every $i\in Q_0$, so in the Tits form \[\ang{\beta,\beta}=\sum_{i\in Q_0}(1-l_i)\beta_i^2-\sum_{i\to j}\beta_i\beta_j\] every monomial appears with a negative coefficient.
		
		We split the proof into several lemmas.	Let $\alpha$ be a Schur root for $Q$.
		
		\begin{lemma}\label{strict}
			\begin{enumerate}[label=(\alph*)]
				\item\label{strict1} We have \[-\ang{\alpha,\alpha}\geq \min_{i\in Q_0}\alpha_i,\] with equality if and only if $Q$ is the $2$-loop quiver and $\alpha=(1)$ or the quiver 	
				\begin{equation}\label{secondcase}
				\begin{tikzcd}
				1 \arrow[r] \arrow[out=45,in=135,loop] & 2 \arrow[out=45,in=135,loop]
				\end{tikzcd}
				\end{equation}\noindent
				and $\alpha=(1,1)$.
				\item\label{strict2} Let $i_0\in Q_0$ satisfy $l_{i_0}\geq 2$, and write $\alpha=\sum_{h=1}^r\beta_h$, for some $\beta_h\in \N^{Q_0}\setminus\set{0}$ and $r\geq 2$. Then
				\[-\sum\ang{\beta_h,\beta_h}\leq -\ang{\alpha,\alpha}-\alpha_{i_0}.\]
			\end{enumerate}
		\end{lemma}
		
		\begin{proof}
			\ref{strict1}. The monomials in the Tits form of $Q$ can only appear with negative coefficient. Since $\alpha_i\alpha_j\geq \alpha_i$ when $\alpha_j\neq 0$, the inequality immediately follows. In order to have equality, it is necessary that the Tits form consists of exactly one monomial. If $l_i\geq 2$ for some $i$, this implies that $Q$ is a $2$-loop quiver, and then it is clear that $\alpha=(1)$ as well. If $l_i=1$ for every $i$, then there are two vertices (just one is excluded, because $Q$ is not the $1$-loop quiver) connected by exactly one arrow, so the quiver is (\ref{secondcase}) and $\alpha=(1,1)$. This proves \ref{strict1}.
			
			\ref{strict2}. If $\alpha_{i_0}\geq 2$, then \[\sum\beta_{h,i_0}^2\leq \alpha_{i_0}^2-2\alpha_{i_0}+2\leq \alpha_{i_0}^2-\alpha_{i_0},\] see \cite[Lemma 6.5]{biswasdhillonhoffmann}.
			Summing this estimate with the trivial inequalities \[\sum\beta_{h,i}\beta_{h,j}\leq \alpha_i\alpha_j,\] one for each arrow $a:i\to j$, gives the conclusion. If $\alpha_{i_0}=1$, then we need only show that \[-\sum\ang{\beta_h,\beta_h}< -\ang{\alpha,\alpha},\] but this is clear because all monomials appear with a positive coefficient and $r\geq 2$.	
		\end{proof}

		\begin{lemma}\label{trdegenoughloops}
			Let $M$ be an indecomposable $\alpha$-dimensional representation over an algebraically closed field that is not a brick. Then \[\trdeg_kk(M)\leq 1-\ang{\alpha,\alpha}-\min_{i\in Q_0}\alpha_i.\]
		\end{lemma}
		
		\begin{proof}
			Using \cite[Corollary 8.4]{scaviaalgebras}, we may write
			\[\trdeg_kk(M)\leq 1-\sum_h\ang{\beta_h,\beta_h}\] where $\beta_h$ is the dimension vector of $\im\phi^{h-1}/\im\phi^{h}$ for a generic $\phi\in\on{End}(M)$. All the entries of $\beta_1$ are non-zero, and since the generic $\phi$ is non-zero there exists a vertex $i_0$ such that $\beta_{2,i_0}\neq 0$. If there is a vertex $i'$ with two loops, then by \Cref{strict}\ref{strict2} \[\trdeg_kk(M)\leq 1-\sum_h\ang{\beta_h,\beta_h}\leq 1-\ang{\alpha,\alpha}-\alpha_{i'}\] and the conclusion follows. Hence we may assume that $l_i=1$ for every $i\in Q_0$. In particular, $Q$ has at least two vertices. If $j\neq i_0$ is another vertex of $Q$, then 
			\begin{align*}
			\alpha_{i_0}\alpha_j=&(\sum_h\beta_{h,i_0})(\sum_{h'}\beta_{h',j})
			\\ =&\sum_h\beta_{h,i_0}\beta_{h,j}+\sum_{h\neq h'}\beta_{h,i_0}\beta_{h',j}
			\\ \geq& \sum_h\beta_{h,i_0}\beta_{h,j}+\beta_{2,i_0}\beta_{1,j}+\beta_{1,i_0}\sum_{h'\geq 2}\beta_{h',j}
			\\ \geq& \sum_h\beta_{h,i_0}\beta_{h,j}+\alpha_j.\end{align*}
			
			Fix an arrow $a:i_0\to i_1$. We consider the estimate above for the term corresponding to $a$ (that is, by letting $j=i_1$), and the inequality
			\[\sum\beta_{h,i}\beta_{h,j}\leq \beta_{i}\beta_{j}\] for every other arrow $i\to j$. Summing up all these inequalities yields
			\[-\sum\ang{\beta_h,\beta_h}\leq -\ang{\alpha,\alpha}-\alpha_j \leq -\ang{\alpha,\alpha}-\min_{i\in Q_0}\alpha_i.\qedhere\]
		\end{proof}

		\begin{lemma}\label{indecnonbrickloops}
			Let $K$ be a field containing $k$. If $M$ is an indecomposable $K$-representation of dimension vector $\alpha$ and is not a brick, then \[\ed_kM\leq -\ang{\alpha,\alpha}.\]
		\end{lemma}
		
		\begin{proof}
			Consider the decomposition $M_{\cl{K}}=\oplus_{h=1}^s N_h$ in indecomposable representations. By \cite[Lemma 12.1]{scaviaalgebras}, this decomposition is defined over $K^{\on{sep}}$, hence over a finite Galois extension $L/K$. Since $M$ is defined over $K$, the Galois group of $L/K$ acts transitively on isomorphism classes of indecomposable summands of $M_L$. We deduce that if some $N_h$ is a brick all the other summands are bricks as well, and that for each $h,h'$ the iterated images of the generic nilpotent endomorphisms of $N_h$ and $N_{h'}$ have the same dimension vectors. We let $\alpha=\dim_KM$, $\beta=\dim_KN_h$, so that $\alpha=s\beta$.
			
			Assume that $N_h$ is a brick for every $h$. Then, since $M$ is not a brick, necessarily $s\geq 2$. We have $\trdeg_kk(N_h)\leq 1-\ang{\beta,\beta}$ by \cite[Corollary 8.4]{scaviaalgebras}. By \Cref{strict}, we have $\min\beta_i\leq-\ang{\beta,\beta}-1$, with the exception of the $2$-loop quiver and $\beta=(1)$, and of the quiver (\ref{secondcase}) and $\beta=(1,1)$. If $\min\beta_i\leq-\ang{\beta,\beta}-1$, using \Cref{estimate} and \cite[Corollary 8.4]{scaviaalgebras}, we obtain: 
			\begin{align*}
			\ed_kM&=\ed_{k(M)}M+\trdeg_kk(M)\\ &\leq\ed_{k(M)}M+\sum_h\trdeg_kk(N_h)\\ &\leq s\min_{i\in Q_0}\beta_i-1+s(1-\ang{\beta,\beta})\\ &\leq -s(1+\ang{\beta,\beta})-1+s(1-\ang{\beta,\beta})\\ &< -2s\ang{\beta,\beta}\leq -s^2\ang{\beta,\beta}=-\ang{\alpha,\alpha}.
			\end{align*}
			If $Q$ is the $2$-loop quiver and $\beta=(1)$, we have $\ang{\beta,\beta}=-1$ and $\ang{\alpha,\alpha}=-s^2$.  If $s\geq 3$, following the same steps as above we obtain
			\[\ed_kM\leq 3s-1<s^2=-\ang{\alpha,\alpha}.\]
			If $s=2$, we may choose a basis so that $M$ is represented by $2$ matrices $A_1,A_2$ commuting with the nilpotent Jordan block of size $2$. This implies that \[A_i=\begin{pmatrix} a_i & 0 \\ b_i & a_i \end{pmatrix},\qquad i=1,2\] so $\ed_kM\leq 4=-\ang{\alpha,\alpha}$. 
			
			If $Q$ is the quiver (\ref{secondcase}) and $\beta=(1,1)$, we have again $\ang{\beta,\beta}=-1$ and $\ang{\alpha,\alpha}=-s^2$. If $s\geq 3$, the same computation yields
			\[\ed_kM\leq 3s-1<s^2=-\ang{\alpha,\alpha}.\]
			If $s=2$, notice that $\phi_{12}:M_1\to M_2$ splits, upon base change to $L$, into the direct sum of two linear maps of the same rank (they are $L$-conjugate), so $\on{rank}\phi_{12}$ is either $0$ or $2$. In the first case $\phi_{12}=0$ and $M$ is the direct sum of two representations of dimension $(2,0)$ and $(0,2)$, and it is easy to see that $\ed_kM\leq 4$. If $\phi_{12}$ is an isomorphism we may identify $M_1$ with $M_2$ via $\phi_{12}$, so that $M$ becomes a representation of the $2$-loop quiver, so $\ed_kM\leq 4$ by the previous case.
			
			Assume now that the $N_h$ are not bricks. Notice that this time $s$ might be $1$. Combining \Cref{estimate} with \Cref{trdegenoughloops}, we get:		 
			\begin{align*}
			\ed_kM&\leq \ed_{k(M)}M+\sum_h\trdeg_kk(N_h)\\ &\leq s\min_{i\in Q_0}\beta_i-1+s(1-\ang{\beta,\beta}-\min_{i\in Q_0}\beta_{i})\\ &<-s\ang{\beta,\beta}+s-1\leq -\ang{\alpha,\alpha},\end{align*}
			the last inequality being equivalent to $-\ang{\beta,\beta}s(s-1)\geq s-1$, which is true because $\alpha$ is sincere and so $\ang{\beta,\beta}=s^{-2}\ang{\alpha,\alpha}<0$. This concludes the proof of \Cref{indecnonbrickloops}.
		\end{proof}
		
		We are ready to prove \Cref{reduction2}. Let $M$ be a $K$-representation that is not a brick, for some field extension $K/k$. If $M$ is indecomposable, then $\ed_kM\leq-\ang{\alpha,\alpha}$  by \Cref{indecnonbrickloops}. If $M$ is decomposable, denote by $M_1,\dots,M_s$ its indecomposable summands, for some $s\geq 2$. By \Cref{estimate} and \cite[Corollary 8.4]{scaviaalgebras} we may write
		\[\ed_kM\leq\ed_{k(M)}M+\sum_h\trdeg_kk(M_h)\leq\min_{i\in Q_0}\alpha_i-\sum_{h=1}^s\ang{\beta_h,\beta_h},\] where $\sum\beta_h=\alpha$.
		To prove that $\ed_kM\leq -\ang{\alpha,\alpha}$, it suffices to show that \[-\ang{\alpha,\alpha}+\sum_h\ang{\beta_h,\beta_h}\geq\min_{i\in Q_0}\alpha_i.\] Assume first that there exists a vertex $j$ such that the sum $\alpha_j=\sum\beta_{h,j}$ has at least two terms, and consider an arrow $a:i\to j$. We have \[\alpha_i\alpha_j-\sum\beta_{h,i}\beta_{h,j}=\sum\beta_{h,i}(\alpha_j-\beta_{h,j})\geq\sum\beta_{h,i}=\alpha_i.\] For every other arrow $a':i'\to j'$, we have \begin{equation}\label{everyother}\alpha_{i'}\alpha_{j'}-\sum\beta_{h,i'}\beta_{h,j'}\geq 0\end{equation} and summing all of these inequalities proves the claim. 
		
		On the other hand, if $\beta_{h,i}\in\set{0,\alpha_i}$ for each vertex $i$ and every $h$, there exist an arrow $a:i\to j$ and two distinct positive integers $h,h'$ such that $\beta_{h,i}=0$ and $\beta_{h',j}=0$. Then \[\alpha_i\alpha_j-\sum\beta_{h,i}\beta_{h,j}=\alpha_i\alpha_j\geq \alpha_i\] and the claim follows by summing this to the inequalities (\ref{everyother}) as in the first case.		
	\end{proof}

	\section{Subquivers}
	
	To finish the proof of \Cref{genericityproperty}, we will need the following combinatorial lemma.
	
	\begin{lemma}\label{combo}
		If $Q$ is not of finite representation type and does not admit at least one loop at every vertex, $Q$ contains a subquiver of one of the following types:
		\begin{enumerate}
			\item\label{quiv1} a tame quiver,
			\item\label{quiv2} a quiver with two vertices and $r\geq 3$ arrows, none of which is a loop and not necessarily pointing in the same direction,
			\[
			\begin{tikzcd}[row sep=1ex]
			& 1 \arrow[r,-] \arrow[r,-,shift left=.6ex, "r"] \arrow[r,-,shift right=.6ex] & 2 
			\end{tikzcd}
			\]
			\item\label{quiv3} a quiver with two vertices, one of which has $s\geq 2$ loops, and with $r\geq 1$ arrows between the two vertices.
			
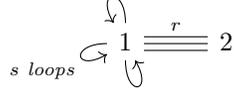
\[
			\begin{tikzcd}
			1 \arrow[out=90,in=120,loop,swap]
			\arrow[out=180,in=210,loop,swap,"s \text{ loops}"]
			\arrow[out=270,in=300,loop,swap]
			\arrow[r,-] \arrow[r,-,shift left=.6ex, "r"] \arrow[r,-,shift right=.6ex] & 2
			\end{tikzcd}
			\]
		\end{enumerate}
	\end{lemma}
	
	\begin{proof}
		Note that $Q$ has at least two vertices, otherwise it would be the trivial quiver with one vertex (which is of finite representation type) or an $r$-loop quiver (which has at least one loop per vertex).
		
		If $Q$ admits at least one loop, we can find two adjacent vertices $i$ and $j$ such that there is at least one loop at $i$ and there are no loops at $j$. If there is exactly one loop at $i$, then $Q$ admits a $1$-loop quiver as a subquiver, and this is of tame representation type. If there are at least two loops at $i$, then $Q$ admits a subquiver of type (\ref{quiv3}). 
		
		Consider the case when that $Q$ does not have any loops. Assume first that there are two vertices $i$ and $j$ connected by $r\geq 2$ arrows. If $r=2$ then $Q$ admits a tame subquiver of type $\widetilde{A}_2$. If $r\geq 3$, then it contains a subquiver of type (\ref{quiv2}). Assume now that $Q$ does not have multiple arrows. If $Q$ admits a cycle, then it admits a tame subquiver of type $\widetilde{A}_n$. It remains to consider the case of a quiver $Q$ without cycles and multiple arrows. By assumption, $Q$ is not of finite representation type. Let $Q'$ be a maximal subquiver of $Q$ that is of finite representation type. Since $Q$ is not of finite representation type, $Q\neq Q'$, and so $Q$ contains a subquiver $Q''$ obtained from $Q'$ by adding one new vertex $j$ to $Q'$, connected only to $i\in Q'_0$ via a unique arrow $i\to j$ (or $j\to i$). One patiently considers all cases for $j$, and concludes that either $Q''$ is of finite representation type, or it contains a tame subquiver. More precisely:
		\begin{itemize}
			\item if $Q'$ is of type $A$, then either $Q''$ is of type $A,D,E$, or it contains a subquiver of type $\widetilde{E}$;
			\item if $Q'$ is of type $D$, then either $Q''$ can be of type $D, E$ or it contains a subquiver of type $\widetilde{D},\widetilde{E}$;
			\item if $Q'$ is of type $E_6$, then $Q''$ either is of type $\widetilde{E}_6,\widetilde{E}_7,E_8$ or contains a subquiver of type $\widetilde{D}_4,\widetilde{D}_5,\widetilde{D}_6$, and similarly in the case when $Q'$ is of type $E_7$ and $E_8$.
		\end{itemize} 
		Assume now that $Q'$ is maximal among subquivers of $Q$ of finite representation type. Then $Q''$ may not be of finite representation type, and so, according to the previous reasoning, it contains a tame subquiver. Therefore, $Q$ contains a tame subquiver.
	\end{proof}
	
	\section{Conclusion of the proof of \Cref{genericityproperty}}	
	Let $Q$ be a quiver. In \Cref{section9} we showed that $R_{Q, \alpha}$ has the genericity property for every dimension vector $\alpha$
	if $Q$ is of finite representation type, or if $Q$ has at least one loop at every vertex. In this section we will establish the converse, thus completing the proof of~\Cref{genericityproperty}. 
	
	Let $Q$ be a quiver such that for every dimension vector $\alpha$, the stack $\mc{R}_{Q,\alpha}$ satisfies the genericity property. Then the same is true for every subquiver of $Q$. This is because a representation of a subquiver $Q'$ may be completed to a representation of the full quiver $Q$ by adding zero spaces and zero linear
	transformations for the vertices and edges in $Q$ but not in $Q'$. This gives rise to
	an isomorphism between the functors $\on{Rep}_{Q, \alpha}$ and $\on{Rep}_{Q', \alpha'}$
	where $\alpha\in \N^{Q_0}$ is obtained from $\alpha'\in\N^{Q_0'}$ by filling in zeros for the missing vertices.
	
	Therefore, it suffices to find for every quiver of the list of \Cref{combo} a dimension vector for which the genericity property does not hold. We will argue in the following way. Suppose that we may find a real root $\alpha$ and a dimension vector $\beta$ such that $\beta_i\leq \alpha_i$ for each vertex $i$ and such that $\ed_k\on{Rep}_{Q,\beta}>0$. By \Cref{deffieldrealroot} we have $\on{ged}_k\mc{R}_{Q,\alpha}=0$, but on the other hand by \cite[Proposition 5.5(b)]{scaviaalgebras} one has $\ed_k\on{Rep}_{Q,\alpha}\geq \ed_k\on{Rep}_{Q,\beta}>0$, so the genericity property does not hold for $\mc{R}_{Q,\alpha}$. 
	
	Consider first the case when $Q$ is a tame quiver, and let $\beta=\delta$ be its null root. By Theorem 7.8(1) of \cite{kirillov} there exists a real root $\alpha$ such that $\alpha_i\geq \delta_i$ for each vertex $i$ of $Q$.
	
	Let now $Q$ be of the second type. The dimension vector $(n,n)$ is a Schur root of generic essential dimension at least $1+(r-1)n^2$, since after fixing an isomorphism between the two vector spaces using one of the arrows, one is reduced to the $(r-1)$-loop quiver. We now construct a suitable real root $\alpha$. One can easily compute the two simple reflections for $Q$:
	\[(x_1,x_2)\mapsto (rx_2-x_1,x_2), \qquad (x_1,x_2)\mapsto (x_1,rx_1-x_2).\]
	If we apply them to $(1,0)$, we get \[(1,0)\mapsto (1,r-1)\mapsto (r^2-r-1,r-1).\]
	Since $r\geq 3$, we have $r^2-r-1>r-1$, hence choosing $\alpha=(r^2-r-1,r-1)$ and $\beta=(r-1,r-1)$ works.
	
	Let now $Q$ be of the third type. Assume first that $r\geq 2$. One can see as in the previous case that the dimension vector $(n,n)$ is a Schur root of generic essential dimension at least $1+(r-1)n^2$. The fundamental region of $Q$ is given by those vectors $(x_1,x_2)$ satisfying \[rx_1-2x_2\geq 0.\]  The vector $(2,1)$ is in the fundamental region and is therefore a Schur root. 	
	There is only one simple reflection, given by \[\sigma:(x_1,x_2)\mapsto (x_1, rx_1-x_2)\]
	From \Cref{genericed} the Schur root $\alpha=(2,2r-1)$ obtained by reflecting $(2,1)$ satisfies: \[\on{ged}_k\mc{R}_{Q,\alpha}=1-\ang{\alpha,\alpha}=2r+4s-4.\]
	On the other hand, since $r\geq 2$, the vector $\beta=(2,2)$ is component-wise smaller than $\alpha$, and \[\on{ged}_k\mc{R}_{Q,\beta}=4r+4s-7>2r+4s-4,\] thus the genericity property does not hold for $\alpha$.
	
	If $r=1$, one may choose $\alpha=(4,3)$ and $\beta=(4,2)$. The vector $\beta$ belongs to the fundamental region $\set{(x_1,x_2): x_1-2x_2\geq 0}$, hence is a Schur root and has generic essential dimension $5$. The vector $\alpha$ is obtained by reflecting $(2,1)$, which belongs to the fundamental region. Hence $\alpha$ is also a Schur root, and has generic essential dimension $4$. This completes the proof of \Cref{genericityproperty}.
	\qed
	
	\begin{example}\label{exloop}	
		Let $r\geq 1$, and consider the $r$-loop quiver $L_r$, here depicted for $r=4$.
		
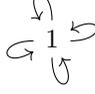
\[
		\begin{tikzcd}
		1 \arrow[out=0,in=30,loop,swap]
		\arrow[out=90,in=120,loop,swap]
		\arrow[out=180,in=210,loop,swap]
		\arrow[out=270,in=300,loop,swap]
		\end{tikzcd}
		\]
		The case $r=1$ has been considered in \Cref{1loop}. Representations of $L_r$ correspond to representations of the free algebra on $r$ generators. It follows from \Cref{genericed} and \Cref{genericityproperty} that \[\ed_k\mc{R}_{L_r,n}\leq 1+(r-1)n^2+\sum_{p}(p^{v_p(n)}-1),\] with equality when \cite[Conjecture 3.10]{biswasdhillonhoffmann} holds for $n$.
		
		This example was originally worked out by Z. Reichstein and A. Vistoli (unpublished). Their proof is in the spirit of \cite{genericity}.
	\end{example}
	
	\section{Proof of \Cref{infiniteged}}
	
	The starting point for the proofs of \Cref{infiniteged} and \Cref{genkron} is the observation that the following inequality holds in the fundamental region of the quiver $Q$.
	
	\begin{lemma}\label{coolestimate}
		Let $Q$ be a quiver, and let $\alpha$ be a dimension vector in the fundamental region of $Q$ such that $\alpha_i>0$ for each vertex $i$. Write $\alpha=\sum_{h=1}^r\beta_h$ for some dimension vectors $\beta_h\in \N^{Q_0}$. Then
		\[-\sum_{h=1}^r\ang{\beta_h,\beta_h}\leq -\ang{\alpha,\alpha}.\] 
		Assume further that for each vertex $i$ there exist at least two vectors $\beta_h$ satisfying $\beta_{hi}\neq 0$. Then \[-\sum_{h=1}^r\ang{\beta_h,\beta_h}\leq -\ang{\alpha,\alpha}-\sum_{i\in Q_0}2(\alpha_i-1)(\sum_{i-j}\frac{\alpha_j}{2\alpha_i}-1).\]
	\end{lemma}
	
	\begin{proof}
		
		Since $\alpha$ belongs to the fundamental region, we have the following inequalities for each vertex $i$: \[(\alpha,e_i)=2\alpha_i-\sum_{i\rightleftarrows j}\alpha_{j}\leq 0.\]
		By algebraic manipulations, starting from:
		\[(\alpha_i\beta_{hj}-\alpha_j\beta_{hi})^2\geq 0\]
		we obtain \[\beta_{hj}\beta_{hi}\leq \frac{\alpha_i}{2\alpha_j}\beta^2_{hj}+\frac{\alpha_j}{2\alpha_i}\beta_{hi}^2\] for every triple $i,j\in Q_0$ and every $h=1,\dots,r$. Hence	
		\begin{align*}
		1-\sum_{h=1}^r\ang{\beta_h,\beta_h}= & 1-\sum_{h,i}\beta_{hi}^2+\sum_{h,i\to j}\beta_{hi}\beta_{hj} \\
		\leq & 1-\sum_{h,i}\beta_{hi}^2+\sum_{h,i\to j}\frac{\alpha_i}{2\alpha_j}\beta_{hj}^2+\sum_{i\to j}\frac{\alpha_j}{2\alpha_i}\beta_{hi}^2\\
		= & 1-\sum_{h,i}\beta_{hi}^2+\sum_{h,i\to j}\frac{\alpha_i}{2\alpha_j}\beta_{hj}^2+\sum_{j\to i}\frac{\alpha_i}{2\alpha_j}\beta_{hj}^2\\
		= & 1+\sum_{h,i}\beta_{hi}^2(\sum_{i-j}\frac{\alpha_j}{2\alpha_i}-1).
		\end{align*}	
		Since $\alpha$ is in the fundamental region, the quantities in the parentheses are non-negative. Since $\sum_{h}\beta_{hi}=\alpha_i$ and $\beta_{hi}\geq 0$, clearly \[\sum_{h}\beta_{hi}^2\leq \alpha_i^2.\] Substituting, we get
		\[1-\sum_{h=1}^r\ang{\alpha_h,\alpha_h}\leq 1+\sum_i\alpha_i^2(\sum_{i- j}\frac{\alpha_j}{2\alpha_i}-1)=1-\ang{\alpha,\alpha}.\]
		Assume now that the conditions of the second part are satisfied. Then one can reach the conclusion using the inequality \[\sum_h\beta^2_{hi}\leq\alpha_i^2-2\alpha_i+2\] for each vertex $i$; see \cite[Lemma 6.5]{biswasdhillonhoffmann}.
	\end{proof}
	
	\begin{lemma}\label{type2}
		Let $r\geq 3$, and $Q$ be a quiver whose underlying graph has the following form:	
		\[
		\begin{tikzcd}[row sep=1ex]
		& 1 \arrow[r,-] \arrow[r,-,shift left=.6ex, "r"] \arrow[r,-,shift right=.6ex] & 2 
		\end{tikzcd}
		\]
		Then, for infinitely many $n$, the genericity property holds for the dimension vector $(n,n)$.
	\end{lemma}
	
	\begin{proof}
		Let $M$ be a $K$-representation of $Q$ of dimension vector $\alpha=(n,n)$, and let $1\leq d\leq 2n$ be the number of indecomposable summands in a Krull-Schmidt decomposition of $M$. By \Cref{estimate} and \cite[Corollary 8.4]{scaviaalgebras}, we may write
		\[\ed_kM\leq n-1+d-\sum\ang{\beta_h,\beta_h}\] for some dimension vectors $\beta_h$ satisfying $\sum\beta_h=\alpha$ (note that the $\beta_h$ are not necessarily the dimension vectors of the summands of $M$).
		By \Cref{coolestimate}, we have
		\[-\sum\ang{\beta_h,\beta_h}\leq (r-2)n^2-4(n-1)(\frac{r}{2}-1)\leq (r-2)n^2-2n+2.\]
		\[\ed_kM\leq d-n+1+(r-2)n^2.\]
		Assume that $n$ is the power of a prime. Then by \Cref{genericed} \[\on{ged}_k\mc{R}_{Q,\alpha}=(r-2)n^2+n.\] If $d\leq 2n-1$, the result follows. On the other hand, if $d=2n$, then $M$ is a direct sum of representations of dimension vectors $(1,0)$ or $(0,1)$, and so $\ed_kM=0$. We conclude that the genericity property holds when $n$ is the power of a prime. 
	\end{proof}
	
	\begin{proof}[Proof of \Cref{infiniteged}]
		By \Cref{genericityproperty}, we may assume that $Q$ does not have at least one loop at every vertex. Moreover, we are allowed to pass to a subquiver of $Q$. By \Cref{combo}, we may assume that $Q$ is of one of the following types:
		
		\begin{enumerate}
			\item\label{qtype1} a quiver obtained from a tame quiver $Q'$ by connecting one extra vertex $i_0$ to exactly one vertex $i_1$ of $Q'$ with $r\geq 1$ arrows,
			\item\label{qtype2} a quiver with two vertices and $r\geq 3$ arrows none of which is a loop and not necessarily pointing in the same direction,
			\item\label{qtype3} a quiver with two vertices, one of which has $s\geq 2$ loops, and with $r\geq 1$ arrows between the two vertices.
		\end{enumerate}  
		Type (\ref{qtype1}) of the list needs further explaination: if the vertex $i_0$ is connected to more than one vertex of $Q'$, it means $Q$ contains a cycle and at least two vertices with at least $3$ edges emanating from them, and so admits a wild proper subquiver $Q''$, and we may consider $Q''$ instead of $Q$ instead.
		
		In case (\ref{qtype3}), $Q$ has a subquiver with at least one loop at every vertex, so the claim holds. Case (\ref{qtype2}) has been treated in \Cref{type2}.
		If $Q$ is of type (\ref{qtype1}), let $\delta$ be the null root of the tame subquiver $Q'$. Fix $m\geq 0$ and define a dimension vector $\alpha$ of $Q$ by setting $\alpha_{i_1}=1$ and $\alpha_i=m\delta_i$ for each $i\neq i_0$. In other words, $\alpha=m\delta+e_{i_0}$, where $\delta$ is viewed as a vector in $\R^{Q_0}$ by extension to zero. Notice that $\alpha$ belongs to the fundamental region of $Q$ for $m\geq 2$, since
		\begin{align*}
		(\alpha,e_i)=\begin{cases}
		2-mr\delta_{i_1} &\text{ if $i=i_0$}\\
		-r &\text{ if $i=i_1$}\\
		0 &\text{ otherwise}.
		\end{cases}
		\end{align*}
		By \cite[Proposition 4.14]{lebruyn} $\alpha$ is an anisotropic Schur root for every $m\geq 2$. We also have
		\[\ang{\alpha,\alpha}=\ang{m\delta,m\delta}+(m\delta,e_{i_0})+\ang{e_{i_0},e_{i_0}}=1-r\alpha_{i_1}.\]
		Therefore, by \Cref{genericed}
		\[\on{ged}_k\mc{R}_{Q,\alpha}=1-\ang{\alpha,\alpha}=r\alpha_{i_1}.\]
		Let now $K$ be a field containing $k$, and let $M$ be an $\alpha$-dimensional $K$-representation of $Q$. By \Cref{estimate}, $\ed_{k(M)}M=0$. We may write $M_{\cl{K}}=M_1\oplus M_2\oplus M_3$, where $M_1$ is the unique indecomposable summand with $(M_1)_{i_0}\neq 0$, $M_2$ is the direct sum of all imaginary indecomposable summands of $M_{\cl{K}}$, and $M_3$ the direct sum of the real ones. Write $\alpha=\beta+c\delta+\gamma$ for the corresponding decomposition of the dimension vector of $M$. By \cite[Corollary 8.4]{scaviaalgebras}, we may write \[\trdeg_kk(M_1)\leq 1-\sum\ang{\beta_h,\beta_h}\] for some decomposition $\beta=\sum\beta_h$.
		Among the $\beta_h$, only one is not supported on the tame subquiver $Q'$, and we denote it by $\beta'$. For every other $\beta_h$, we have $\ang{\beta_h,\beta_h}=0$. Writing $\beta'=e_{i_0}+\beta''$, for some $\beta''\in\R^{Q'_0}$, we obtain \[\trdeg_kk(M_1)\leq 1-\ang{\beta',\beta'}=1-\ang{e_{i_0},e_{i_0}}-\ang{\beta'',\beta''}-(\beta'',e_{i_0})\leq 1+r\beta'_{i_1}.\]
		From \Cref{tamegen}, $\trdeg_kk(M_2)\leq c$ and by \cite[Remark 9.1]{scaviaalgebras}, $\trdeg_kk(M_3)=0$. Thus
		\[\ed_kM=\trdeg_kk(M)\leq r\beta'_{i_1}+c\leq r(\beta_{i_1}+c\delta_{i_1})\leq r\alpha_{i_1}.\]
		Therefore, the genericity property holds for the dimension vector $\alpha$.
	\end{proof}
	
	\begin{rmk}\label{moreprecise}
		If $Q$ is of type (\ref{qtype1}) in the list above, $\delta$ is the null root of its tame subquiver $Q'$, $\alpha=m\delta+e_{i_0}$, we have just shown that (when $m\geq 2$) $\alpha$ is a Schur root and the genericity property holds for $\alpha$. If $Q$ is of type (\ref{qtype2}), we have shown in the proof of \Cref{type2} that the genericity property holds for $(n,n)$ when $n$ is the power of a prime. Finally, if $Q$ is of type (\ref{qtype3}), it contains the $s$-loop quiver for $s\geq 2$ as a subquiver with unique vertex $i_0$, and so by \Cref{genericityproperty} the genericity property holds for $me_{i_0}$ for every $m\geq 0$.
		
		By \Cref{combo}, every wild quiver contains at least one subquiver of type (\ref{qtype1}), (\ref{qtype2}) or (\ref{qtype3}). To produce Schur roots for which the genericity property holds, it thus suffices to identify one of these subquivers. 
	\end{rmk}

	\section{Proof of \Cref{genkron}}
	
	For a positive integer $r$, let $K_r$ be the $r$-Kronecker quiver.
	
	Let $\alpha=(a,b)$ be a dimension vector for $K_r$. The quiver $K_1$ is of finite representation type, so \[\ed_k\on{Rep}_{K_1,\alpha}=0.\] The indecomposable representations of $K_2$ were classified by Kronecker (see \cite[Theorem 3.6]{burban2012two} for a description over an arbitrary field). It follows from the classification that \[\ed_k\on{Rep}_{K_2,\alpha}=\floor{\frac{a+b}{2}}.\]
	
	The main result of this section is the proof \Cref{genkron}. Recall that we have already shown in the course of proving \Cref{genericityproperty} that the genericity property fails for the Schur root $(r^2-r-1,r-1)$. Therefore one cannot expect \Cref{genkron} to hold for each Schur root.
	
	\begin{proof}[Proof of \Cref{genkron}]
		The argument follows steps similar to those of the proof of \Cref{reduction2}.

		\begin{lemma}\label{kronalg}
			Assume that $M$ is an indecomposable $\alpha$-dimensional representation of $K_r$ over an algebraically closed field $K$, and that $M$ is not a brick. Then \[\trdeg_kk(M)\leq a^2+b^2-rab-\min(a,b).\]	
		\end{lemma}
		\begin{proof}
			Let $\phi\in\on{End}M$ be a generic nilpotent endomorphism of $M$. Write \[\alpha_h=(a_h,b_h):=\dim\on{Im}\phi^{h-1}/\dim\on{Im}\phi^h\] for every $h\geq 0$. If $a_1=a$, this means that there exists a nilpotent endomorphism $\psi$ of $M$ such that $\psi_1=0$ and $\psi_2\neq 0$. We may choose bases of $M_1$ and $M_2$ in such a way that $\psi_2$ is represented by a nilpotent matrix in Jordan form. With respect to these bases, the matrices $A_1,\dots,A_r$ corresponding to the $r$ arrows of $K_r$ all have at least one common row made of only zeros. This is impossible, since $M$ was supposed to be indecomposable. An analogous reasoning proves that $b_1\neq b$, so each of the decompositions $a=\sum a_h$ and $b=\sum b_h$ contains at least two summands. Using \cite[Corollary 8.4]{scaviaalgebras} and \Cref{coolestimate}, we obtain	
			\[\trdeg_kk(M)\leq 1-\sum_{h=1}^r\ang{\alpha_h,\alpha_h}\leq 1-\ang{\alpha,\alpha}-f(a,b)\]
			where \[f(a,b)=2(a-1)(\frac{rb}{2a}-1)+2(b-1)(\frac{ra}{2b}-1).\]
			By \Cref{effe} below, $f(a,b)\geq \min(a,b)-1$. Therefore \[\trdeg_kk(M)\leq 1-a^2-b^2+rab-\min(a,b)+1.\]
			On the other hand, $\ed_{k(M)}M\leq \min(a,b)-1$ by \Cref{estimate}, hence \[\ed_kM\leq 1-\ang{\alpha,\alpha}.\qedhere\]
		\end{proof}
		
		\begin{lemma}\label{effe}
			Assume that $\alpha=(a,b)$ is in the fundamental region of $K_r$, $r\geq 3$. Let \[f(a,b)=2(a-1)(\frac{rb}{2a}-1)+2(b-1)(\frac{ra}{2b}-1).\] Then \[f(a,b)\geq\min(a,b)-1.\]
		\end{lemma}
		
		\begin{proof}
			Since $(a,b)$ belongs to the fundamental region of $K_r$, we have $2a\leq rb$ and $2b\leq ra$. 
			Moreover, since $f$ is symmetric, we may assume that $a\geq b$. Then $\frac{ra}{2b}\geq \frac{r}{2}$, so \[f(a,b)\geq 2(b-1)(\frac{ra}{2b}-1)\geq (b-1)(r-2)\geq b-1.\qedhere\]
		\end{proof}
		
		\begin{lemma}\label{indecnonbrickkron}
			Assume that $M$ is an indecomposable $\alpha$-dimensional representation of $K_r$ over an arbitrary field $K$ containing $k$. If $M$ is not a brick, then \[\ed_kM\leq 1-\ang{\alpha,\alpha}.\]
		\end{lemma}
		\begin{proof}
			Consider the decomposition $M_{\cl{K}}=\oplus_{h=1}^s N_h$ in indecomposable representations. By \cite[Lemma 12.1]{scaviaalgebras}, this decomposition is defined over $K^{\on{sep}}$, hence over a finite Galois extension $L/K$. Since $M$ is indecomposable, the Galois group of $L/K$ acts transitively on isomorphism classes of indecomposable summands of $M_L$. We deduce that if one of the $N_h$ is a brick all of them are, and that for each $h,h'$ the iterated images of the generic nilpotent endomorphisms of $N_h$ and $N_{h'}$ have the same dimension vectors. We let $\alpha=\dim_KM$, $\beta=(\beta_1,\beta_2)=\dim_KN_h$, so that $\alpha=s\beta$.
			
			Assume that $N_h$ is a brick for every $h$. Then, since by assumption $M$ is not a brick, necessarily $s\geq 2$. We have $\trdeg_kk(N_h)\leq 1-\ang{\beta,\beta}$ by \cite[Corollary 8.4]{scaviaalgebras}. Since $\beta$ is in the fundamental region of $K_r$, it satisfies the inequalities \[2\beta_1-r\beta_2\leq 0, \qquad 2\beta_2-r\beta_1\leq 0,\] which imply \[-\ang{\beta,\beta}=-\beta_1^2-\beta_2^2+r\beta_1\beta_2\geq\max(\beta_1^2-\beta_2^2,\beta_2^2-\beta_1^2).\] If $\beta_1\neq \beta_1$, we obtain $-\ang{\beta,\beta}\geq \min\beta_i$, which is also true if $\beta_1=\beta_2$. We use \cite[Corollary 8.4]{scaviaalgebras} to obtain: 
			\begin{align*}
			\ed_kM&=\ed_{k(M)}M+\trdeg_kk(M)\\ &\leq\ed_{k(M)}M+\sum_h\trdeg_kk(N_h)\\ &\leq s\min_{i\in Q_0}\beta_i-1+s(1-\ang{\beta,\beta})\\ &\leq -2s\ang{\beta,\beta}+s-1\\ &\leq 1-s^2\ang{\beta,\beta}=1-\ang{\alpha,\alpha}.
			\end{align*}
			The last inequality holds because it is equivalent to $-s(s-2)\ang{\beta,\beta}\geq s-2$.
			Assume that the $N_h$ are not bricks. We still have $-\ang{\beta,\beta}\geq\min\beta_i$, this time using \Cref{coolestimate} instead of \cite[Corollary 8.4]{scaviaalgebras} to prove $\trdeg_kk(N_h)\leq 1-\ang{\beta,\beta}$. Notice that this time $s$ might be $1$. Combining \Cref{kronalg} with \Cref{estimate}, we get:		 
			\begin{align*}
			\ed_kM&\leq \ed_{k(M)}M+\sum_h\trdeg_kk(N_h)\\ &\leq s\min_{i\in Q_0}\beta_i-1+s(1-\min_{i\in Q_0}\beta_{i}-\ang{\beta,\beta})\\ &\leq-s^s\ang{\beta,\beta}= -\ang{\alpha,\alpha},\end{align*}
			the last inequality being equivalent to $-s(s-1)\ang{\beta,\beta}\geq s-1$.
		\end{proof}
		Let $K$ be a field extension of $k$, and let $M$ be an $\alpha$-representation of $K_r$ that is not a brick. If $M$ is indecomposable, $\ed_kM\leq 1-\ang{\alpha,\alpha}$ by \Cref{indecnonbrickkron}. If $M$ is decomposable, set $\alpha=(a,b)$, and write $M=\oplus_{h=1}^sM_h$ for the decomposition of $M$ in indecomposable representations, where $s\geq 2$. Let $\beta_h=(a_h,b_h)=\dim M_h$. If there are $h,h'$ such that $\beta_h=(a_h,0)$ and $\beta_{h'}=(0,b_{h'})$, then \[\ang{\beta_h+\beta_{h'},\beta_h+\beta_{h'}}\leq \ang{\beta_h,\beta_h}+\ang{\beta_{h'},\beta_{h'}},\] so we reduce to smaller $a$ and $b$. Now say that there are no $\beta_h$ of the form $(0,b_h)$ (the other case is symmetric). If those of the form $(a_h,0)$ sum to $\gamma=(c,0)$, we see that \[\ang{\alpha-\gamma,\alpha-\gamma}\geq\ang{\alpha,\alpha}\] since this reduces to $c(2a-rb-c)\leq 0$, which is true because $2a-rb\leq 0$. Therefore, we may assume that each $M_h$ has dimension vector with both entries different from zero for each $h$. In this case, by \Cref{estimate}, \cite[Corollary 8.4]{scaviaalgebras}, \Cref{coolestimate} and \Cref{effe}, we obtain:
		\begin{align*}\ed_kM\leq&\sum_h\ed_kM_h \\ =&\sum_{h}(\ed_{k(M_h)}M_h+\trdeg_kk(M_h))\\ \leq&\sum_{h}(\min(a_h,b_h)-\ang{\beta_h,\beta_h})\\ \leq& \min(a,b)-\ang{\alpha,\alpha}-f(\alpha)\\ \leq& 1-\ang{\alpha,\alpha}.\qedhere\end{align*} 
	\end{proof}
	
	\section*{Acknowledgements}
	I am very grateful to Angelo Vistoli for proposing to me this topic of research, and for making me aware of the methods in \cite{biswasdhillonhoffmann}, and to my advisor Zinovy Reichstein for his help and suggestions. I thank Roberto Pirisi and Mattia Talpo for helpful comments, and Ajneet Dhillon and Norbert Hoffmann for helpful correspondence.

	\bibliography{bibliografia}

\begin{thebibliography}{10}

\bibitem{bensonreichstein}
Dave Benson and Zinovy Reichstein.
\newblock Fields of definition for representations of associative algebras.
\newblock {\em arXiv:1702.06447, to appear in Proceedings of the Edinburgh
  Math. Society}.

\bibitem{berhuy2005notion}
Gr{\'e}gory Berhuy and Zinovy Reichstein.
\newblock On the notion of canonical dimension for algebraic groups.
\newblock {\em Advances in Mathematics}, 198(1):128--171, 2005.

\bibitem{biswasdhillonhoffmann}
Indranil Biswas, Ajneet Dhillon, and Norbert Hoffmann.
\newblock On the essential dimension of coherent sheaves.
\newblock {\em J. Reine Angew. Math.}, 735:265--285, 2018.

\bibitem{bourbakialgebra}
Nicolas Bourbaki.
\newblock {\em Algebra}.
\newblock Springer, 1998.

\bibitem{genericity0}
Patrick Brosnan, Zinovy Reichstein, and Angelo Vistoli.
\newblock Essential dimension of moduli of curves and other algebraic stacks.
  {W}ith an appendix by {N}. {F}akhruddin.
\newblock {\em J. European Math. Society (JEMS) 13 (2011), no. 4, 1079-1112}.

\bibitem{burban2012two}
Igor Burban and Olivier Schiffmann.
\newblock Two descriptions of the quantum affine algebra
  {$U_v(\hat{\mathfrak{sl}}_2)$} via hall algebra approach.
\newblock {\em Glasgow Mathematical Journal}, 54(2):283--307, 2012.

\bibitem{cernele2015essential}
Shane Cernele and Zinovy Reichstein.
\newblock Essential dimension and error-correcting codes.
\newblock {\em Pacific Journal of Mathematics}, 279(1):155--179, 2015.

\bibitem{colliothelenekarpenkomerkurjev}
Jean-Louis Colliot-Th{\'e}l{\`e}ne, Nikita Karpenko, and Alexander Merkurjev.
\newblock Rational surfaces and the canonical dimension of {PGL}6.
\newblock {\em St. Petersburg Mathematical Journal}, 19(5):793--804, 2008.

\bibitem{curtis1966representation}
C.W. Curtis and I.~Reiner.
\newblock {\em Representation Theory of Finite Groups and Associative
  Algebras}.
\newblock AMS Chelsea Publishing Series. Interscience, 1966.

\bibitem{drozd}
Yuriy Drozd.
\newblock Tame and wild matrix problems.
\newblock {\em Representations and quadratic forms (Russian), Akad. Nauk
  Ukrain. SSR, Inst. Mat., Kiev, 1979, 39-74, 154.}

\bibitem{hoffmann}
Norbert Hoffmann.
\newblock Rationality and {P}oincare families for vector bundles with extra
  structure on a curve.
\newblock {\em arXiv:math/0511656}, 2005.

\bibitem{kac1}
Victor~G. Kac.
\newblock Infinite root systems, representations of graphs and invariant theory
  i.
\newblock {\em Inventiones mathematicae 56 (1980): 57-92}, 1980.

\bibitem{kac2}
Victor~G. Kac.
\newblock Infinite root systems, representations of graphs and invariant theory
  ii.
\newblock {\em Journal of Algebra 78, 141-162}, 1982.

\bibitem{karpenko2013upper}
Nikita~A Karpenko.
\newblock Upper motives of algebraic groups and incompressibility of
  severi--brauer varieties.
\newblock {\em Journal f{\"u}r die reine und angewandte Mathematik (Crelles
  Journal)}, 2013(677):179--198, 2013.

\bibitem{karpenko2006canonical}
Nikita~A Karpenko and Alexander~S Merkurjev.
\newblock Canonical p-dimension of algebraic groups.
\newblock {\em Advances in Mathematics}, 205(2):410--433, 2006.

\bibitem{karpenko2014numerical}
Nikita~A. Karpenko and Zinovy Reichstein.
\newblock A numerical invariant for linear representations of finite groups,
  with an appendix by {J}ulia {P}evtsova and {Z}inovy {R}eichstein.
\newblock {\em Commentarii Math. Helveciti, Volume 90, no. 3, 2015, 667-701.}

\bibitem{king}
Alastair~D. King.
\newblock Moduli of representations of finite dimensional algebras.
\newblock {\em The Quarterly Journal of Mathematics 45.4: 515-530.}, 1994.

\bibitem{kirillov}
Alexander Kirillov.
\newblock {\em Quiver Representations and Quiver Varieties}.
\newblock Graduate Studies in Mathematics. American Mathematical Society, 2016.

\bibitem{laumonmoretbailly}
G\'erard Laumon and Laurent Moret-Bailly.
\newblock {\em Champs alg\'ebriques}.
\newblock Springer, 2000.

\bibitem{lebruyn}
Lieven LeBruyn.
\newblock {\em Noncommutative Geometry and Cayley-smooth Orders}.
\newblock Pure and Applied Mathematics 290, Chapman and Hall, 2007.

\bibitem{merkurjev2013essential}
Alexander~S Merkurjev.
\newblock Essential dimension: a survey.
\newblock {\em Transformation groups}, 18(2):415--481, 2013.

\bibitem{reichsteinwhatis}
Zinovy Reichstein.
\newblock What is... essential dimension?
\newblock {\em Notices of the American Mathematical Society}, 2012.

\bibitem{genericity}
Zinovy Reichstein and Angelo Vistoli.
\newblock A genericity theorem for algebraic stacks and essential dimension of
  algebraic hypersurfaces.
\newblock {\em Journal of the European Math. Society (JEMS), 15 (2013), no. 6,
  1999-2026.}

\bibitem{scaviaalgebras}
Federico Scavia.
\newblock New invariants for representations of algebras.
\newblock {\em arXiv:1804.00642 [math.AG]}.

\bibitem{schiffler}
Ralf Schiffler.
\newblock {\em Quiver Representations}.
\newblock CMS Books in Mathematics. Springer, 2006.

\bibitem{deffieldrealroot}
Aidan Schofield.
\newblock The field of definition of a real representation of a quiver {Q}.
\newblock {\em Proceedings of the American Mathematical Society},
  116(2):293--295, 1992.

\bibitem{schofield91}
Aidan Schofield.
\newblock General representations of quivers.
\newblock {\em Proceedings of the London Mathematical Society},
  s3-65(1):46--64, 1992.

\bibitem{serre2012linear}
Jean-Pierre Serre.
\newblock {\em Linear representations of finite groups}, volume~42.
\newblock Springer Science \& Business Media, 2012.

\end{thebibliography}
	\bibliographystyle{plain}
	
\end{document}